\documentclass[article,11pt]{amsart}

\usepackage{amsmath,amsthm,amssymb,amscd}

\usepackage[all,cmtip,color]{xy}
\usepackage[german,english]{babel}
\usepackage{amsmath}

\usepackage{amssymb}

\usepackage{xcolor}
\usepackage{hyperref}

\hypersetup{
  colorlinks,
  citecolor=blue,
  linkcolor=blue,
  urlcolor=blue}
  
\usepackage{mathrsfs} 
\usepackage{bbm}
\usepackage{color}
\usepackage{dsfont}
\usepackage{todonotes}
\usepackage{tikz}
\usetikzlibrary{cd,arrows,matrix,positioning}

\usepackage[margin = 2.8cm,headsep=1cm]{geometry}

\usepackage[all,cmtip]{xy}
\usepackage{colonequals,enumerate}
\usepackage{enumitem} 
\usepackage[normalem]{ulem}

 \usepackage{comment}

\def\sheafhom{\mathcal{H}om}

\DeclareMathOperator{\Ext}{\mathsf{Ext}}

\DeclareMathOperator{\A}{\mathsf{A}}
\DeclareMathOperator{\B}{\mathsf{B}}

\renewcommand{\lim}{\mathsf{lim}}

\DeclareMathOperator{\Aut}{\mathsf{Aut}}

\DeclareMathOperator{\G}{\mathbb{G}}

\DeclareMathOperator{\GL}{GL}

\DeclareMathOperator{\M}{M}

\DeclareMathOperator{\R}{R}

\DeclareMathOperator{\Spec}{\mathsf{Spec}}

\newcommand{\BA}{{\mathbb{A}}}

\newcommand{\BC}{{\mathbb{C}}}

\newcommand{\BF}{{\mathbb{F}}}
\newcommand{\BG}{{\mathbb{G}}}

\newcommand{\BM}{{\mathbb{M}}}

\newcommand{\BP}{{\mathbb{P}}}
\newcommand{\BQ}{{\mathbb{Q}}}

\newcommand{\BZ}{{\mathbb{Z}}}

\newcommand{\FC}{{\mathcal C}}

\newcommand{\FM}{{\mathcal M}}

\newcommand{\FO}{{\mathcal O}}

\newcommand{\FU}{{\mathcal U}}

\newcommand{\Pic}{\mathrm{Pic}}
\newcommand{\PIC}{\mathfrak{Pic}}

\renewcommand{\M}{\mathrm{M}}
\renewcommand{\B}{\mathrm{B}}

\newcommand{\pBPS}{p\mathrm{BPS}}
\newcommand{\BPS}{\mathrm{BPS}}
\newcommand{\DCGS}{\mathrm{DCGS}}

\newtheorem{theorem}[equation]{Theorem}
\newtheorem{proposition}[equation]{Proposition}

\newtheorem{corollary}[equation]{Corollary}

\newtheorem{lemma}[equation]{Lemma}
\theoremstyle{definition}
\newtheorem{definition}[equation]{Definition}

\newtheorem{rmk}[equation]{Remark}

\newtheorem{example}[equation]{Example} 
\newtheorem{assumption}[equation]{Assumption}

\newtheorem{setting}[equation]{Setting}

\numberwithin{equation}{subsection}

\begin{document}

\author{Francesca Carocci, Giulio Orecchia, Dimitri Wyss}
\email{francesca.carocci@unige.ch}
\address{University of Geneva, Mathematics Department, Office 6.05,
Rue du Conseil-Général 7-9  CH-1204 Genève, Switzerland}

\email{giulioorecchia@gmail.com }
\address{White Oak Asset Management,  Rue du Rhône 17, CH-1204, Geneva}

\email{dimitri.wyss@epfl.ch}
\address{Ecole Polytechnique F\'ed\'erale de Lausanne (EPFL), Chair of Arithmetic Geometry, CH-1015 Lausanne, Switzerland
}

\title{BPS invariants from $p$-adic integrals}

\keywords{Moduli of sheaves, Donaldson--Thomas invariants, p-adic integration}
\begin{abstract} 

We define  $p$-adic $\BPS$ or $\pBPS$-invariants for moduli spaces $\M_{\beta,\chi}$ of 1-dimensional sheaves on del Pezzo and K3 surfaces by means of integration over a non-archimedean local field $F$.
Our definition relies on a canonical measure $\mu_{can}$ on the $F$-analytic manifold associated to $\M_{\beta,\chi}$ and the $\pBPS$-invariants are integrals of natural $\BG_m$-gerbes with respect to $\mu_{can}$. A similar construction can be done for meromorphic and usual Higgs bundles on a curve.

Our main theorem is a $\chi$-independence result for these $\pBPS$-invariants. For 1-dimensional sheaves on del Pezzo surfaces and meromorphic Higg bundles, we obtain as a corollary the agreement of $\pBPS$ with usual $\BPS$-invariants trough a result of Maulik–Shen \cite{maulik2020cohomological}.

\end{abstract}
\maketitle

\section{Introduction}

Donaldson--Thomas (DT)--invariants, first introduced in \cite{thomasholomorphic}, count stable sheaves with some fixed Chern character $\gamma\in H^*(X,\mathbb Z)$ on a smooth Calabi--Yau  3-fold $X$.  While with the original machinery the invariants could only be defined for moduli of sheaves where no strictly semi-stables occur, Joyce--Song \cite{joycetheory} and Kontsevich--Soibelman \cite{kontsevichsoibelman} independently developed a generalised theory allowing the definition of  DT--invariants for moduli of objects (in CY3 categories) also in those cases where strictly semi-stables are present. 
The two  approaches have some differences and the resulting generalised numerical invariants, denoted in the above references by $\bar{\text{DT}}_{\gamma}$  and  $\hat{\text{DT}}_{\gamma}$ or $\Omega(\gamma)$ respectively do not coincide. The relation between the two is however understood and explained for example in \cite[\S~6.2]{joycetheory} or \cite[\S~6.7]{davisondimone}.

Both theories admit, at least in some specific geometries, 
refinements to motivic, cohomological and sheaf theoretic invariants. 
Let $\BM_{\gamma}(X)$ denote some moduli stack of semi-stable sheaves on a smooth Calabi--Yau  3-fold; Joyce--Song theory is refined by the cohomology of a certain perverse sheaf  $\mathcal D T _{\gamma}$ on the moduli stack  obtained (given an orientation data) gluing locally defined vanishing cycle sheaves (see \cite{joycedarboux}, \cite{brav2015symmetries}). 
 Kontsevich--Soibelman theory should instead be refined by the cohomology $\text{H}^*(\M_{\gamma}(X),\Phi_{\gamma})$ for  $\Phi_{\gamma}$ the so-called \emph{BPS-sheaf}, named after Bogomol’nyi-Prasad-Sommerfield, defined on the moduli space $\M_{\gamma}(X)$ of S-equivalence classes. For $X$ a compact CY3 threefold
the existence of the sheaf $\Phi_{\gamma}$  is still conjectural; see \cite{todagopakumar} for a conjectural definition inspired to the case of moduli spaces of quiver representation with potential. In the latter situation,  the BPS sheaf  $\Phi_{\gamma}$, the sheaf on the moduli stack  $\mathcal D T _{\gamma}$  as well their relation are understood  \cite{davisonintegrality}.

Following \cite{joycetheory} and \cite{davisonintegrality} in this paper we call BPS--invariants the generalised DT invariants of Kontsevich--Soibelman, since as suggested by the authors themselves, their invariants should count BPS states.

Particularly interesting from the perspective of enumerative geometry is the case of one dimensional sheaves on $X$, i.e. sheaves with Chern character $\gamma=(0,0,\beta,\chi).$ In this case the invariants arising from the moduli spaces of Gieseker semi-stable sheaves are (conjecturally) related to Gopakumar--Vafa invariants \cite{katz},\cite{maulikgopakumar},\cite{todagopakumar}. 

 Our main focus is the special case where the CY3-fold is a del Pezzo or a K3 local surface
  $X=\operatorname{Tot}(S)$. In these cases the BPS sheaf is better understood. For a local del Pezzo case a result of Meinhardt \cite{meinhomone} states that  $\Phi_{\beta,\chi}$ exists and coincides with the intersection complex $IC_{\M_{\beta,\chi}}$ of the moduli space of S-equivalence classes. For a local K3 the situation is more complicated and an explicit description of the BPS sheaf has only recently been given in \cite{davison2023bps}.

The main result of the present paper gives a surprising relation between the $\BPS$-invariants of local del Pezzo surfaces and certain non-archimedean (or $p$-adic) integrals. The relation is indirect and relies on results of Maulik--Shen \cite{maulik2020cohomological} which we explain next.


\subsection{Cohomological $\chi$--independence and trace of Frobenius}
Let $S$ be a smooth  surface over $\mathbb C$ and $\beta$ an ample base point free class. We denote by   $\BM_{\beta,\chi}(S)$ the moduli stack of Gieseker semi-stable (with respect to some fixed polarization $H$) $1$-dimensional sheaves, and by $\M_{\beta,\chi}(S)$  the moduli space of S-equivalence classes.
When $S$ is del Pezzo, the stack of semi-stables is smooth. Under this hypothesis Maulik--Shen \cite{maulik2020cohomological} have recently proved the independence of the intersection cohomology from the Euler characteristic:\footnote{Strictly speaking in \cite{maulik2020cohomological} $S$ is assumed to be toric, but with the recent results of \cite{yuan} this assumption can be dropped, see \cite[Remark 0.2]{maulik2020cohomological} and Corollary \ref{cor:irr} below.}
\[\operatorname{IH}^*(\M_{\beta,\chi}(S))\cong \operatorname{IH}^*(\M_{\beta,\chi'}(S))\;\;\forall\;\chi,\chi'.\]
In fact they prove a stronger statement, namely that the pushforward of the intersection complexes along the Hilbert-Chow morphisms $h_{\chi},h_{\chi'}$ are isomorphic. We denoted by $h_\chi$ the morphism
\[ h_{\chi}:\M_{\beta,\chi}(S) \rightarrow \B = \BP H^0(S,\FO_S(\beta)),\]
associating to a  sheaf its Fitting support.

The result  proves a refinement of a special case of a more general and long standing conjecture in the enumerative geometry of smooth CY3-fold,
known as the Pandharipande--Thomas strong rationality conjecture \cite{PT} and later reformulated by Toda \cite{todamultiple} as the multiple cover formula conjecture for generalised DT--invariants of moduli spaces of one dimensional sheaves. Translated in the Kontsevich--Soibelmain theory, the conjecture predicts the independence of the 
BPS--invariants from the Euler characteristic.

 The conjecture is expected to hold more in general at the refined level.



 

Now, choosing a spreading out, we can assume that the moduli spaces $\M_{\beta,\chi}(S)$ we are considering are defined over some large finite field $k=\mathbb F_q$ and that $\chi$--independence of $\Phi_{\beta,\chi}$ holds over $k$ (see for example \cite[Section 6]{BBDG18} for precise definitions and results on spreading out of constructible complexes). We can then look at the function:
 \[\operatorname{BPS}_{\beta,\chi}\colon \M_{\beta,\chi}(k)\to\mathbb C,\qquad x\mapsto q^{-\dim \FM_{\beta,\chi}} \operatorname{Tr}(\operatorname{Fr},\Phi_{\beta,\chi, x}).\]

We call a pair $(\beta,\chi)$ generic (with respect to $H$), if any Gieseker semi-stable sheaf in $\BM_{\beta,\chi}$ is stable. In this case $\BM_{\beta,\chi} \to \M_{\beta,\chi}$ is a $\BG_m$-gerbe and in particular $\M_{\beta,\chi}$ is smooth. Then the function $\operatorname{BPS}_{\beta,\chi}$ will have the following two properties:
 \begin{enumerate}
 \item$\BPS_{\beta,\chi}\equiv q^{-\dim \M_{\beta,\chi}} $ if $(\beta,\chi)$ is generic;
 \item For all $\chi,\chi'\in\mathbb Z$ and for all $y\in B(k)$ we have 
 \[\sum _{x\in h_{\chi}^{-1}(y)(k)} \operatorname{BPS}_{\beta,\chi}(x)=\sum _{x\in h_{\chi'}^{-1}(y)(k)} \operatorname{BPS}_{\beta,\chi'}(x),\] 
   \end{enumerate}
the second one being a consequence of \cite{maulik2020cohomological}.

The main content of the paper is the construction of a function $\pBPS_{\beta,\chi}$ via non-archimedean integration satisfying these two properties. 

\subsection{p-adic BPS function and its invariance}
Consider $S$ a smooth  del Pezzo surface over $\Spec(\FO)$, where $\FO$ denotes the ring of integers of a non-archimedean local field $F$ with residue field $k \cong \BF_q$.

As recalled above, the del Pezzo hypothesis ensures that the moduli stack is smooth, which allows us to construct in  Section~\ref{pbasics} a \emph{canonical measure} $\mu_{can}$ on the $F$-analytic manifold $\M_{\beta,\chi}(\mathcal O)^\natural$ associated to $\M_{\beta,\chi}$. 

In fact, the construction of $\mu_{can}$ works more generally for normal, generically stabiliser-free, Artin stacks $\mathcal M$ admitting a universally closed morphism $\mathcal M\xrightarrow{\pi}\M$ to a quasi-projective variety $\M$ such that $\pi$ is generically an isomorphism. The existence of a canonical measure in this context should be of independent interest.

Importantly, these hypothesis are also satisfied by moduli stacks $\M_{\beta,\chi}$ of semi-stable sheaves on $S\to\Spec(\FO)$ a smooth projective K3 over $\FO$ for many choices of $(\beta,\chi)$ (see \S~3.2 for details), as well as by moduli stacks of (usual) Higgs bundles.

Once defined the canonical measure on $\M_{\beta,\chi}(\mathcal O)^\natural$,
we then define the non-archimedean BPS-function $\pBPS_{\beta,\chi}$ as follows:
for any $x \in \M_{\beta,\chi}(k)$ denote by $\M_{\beta,\chi}(\FO)_x \subset M_{\beta,\chi}(\FO)^\natural$ the ball of $\FO$-rational points specialising to $x$ over $k$. Then 
\[\pBPS_{\beta,\chi}:\M_{\beta,\chi}(k) \rightarrow \BC\]
is given by:
\begin{equation}\label{pbpsdef}
\pBPS_{\beta,\chi}(x) =q^{-\dim M_{\beta,\chi}}  \int_{M_{\beta,\chi}(\FO)_x} \varphi_{\beta,\chi}^g d\mu_{can},\end{equation}
where $\varphi_{\beta,\chi}$ is a certain complex valued function associated with the natural $\BG_m$-gerbe coming from the $\BG_m$-rigidification of $\BM_{\beta,\chi}$ (see Section~\ref{sec:mainsec} for the precise definition) and $g$ is the arithmetic genus of the curves in the linear system $|\beta|$.

Our main result says that the non-archimedean BPS function enjoys the two properties above:
\begin{theorem}[\ref{mainthm}]\label{mainthmintro}  Let  $S \to \Spec(\FO)$ be either a smooth projective del Pezzo surface or a K3 surface satisfying Assumption \ref{K3as}. The function $\pBPS_{\beta,\chi}: \M_{\beta,\chi}(k) \rightarrow \BC$ satisfies the following two properties:
 \begin{enumerate}
 \item$\pBPS_{\beta,\chi} \equiv q^{-\dim \M_{\beta,\chi}} $ if $(\beta,\chi)$ is generic;
 \item For all $\chi,\chi'\in\mathbb Z$ and for all $y\in B(k)$ we have 
 \[\sum _{x\in h_{\chi}^{-1}(y)(k)} \pBPS_{\beta,\chi}(x)=\sum _{x\in h_{\chi'}^{-1}(y)(k)} \pBPS_{\beta,\chi'}(x).\]
\end{enumerate}
\end{theorem}
\subsection{Relation to previous works, consequences and final considerations }

The original motivation for Theorem~\ref{mainthmintro} is a degree-independence result for Higgs bundles of rank coprime to the degree \cite[Theorem 7.15]{GWZ20}, which in turn is a special case of a conjecture by Mozgovoy-Schiffmann \cite{mozgovoy2014counting}\footnote{There is an alternative proof by Yu  \cite{Yu18} and proofs of the complete conjecture by Mellit \cite{Me20} and Kinjo-Koseki \cite{KK21}.}. Thus not surprisingly, also the proof of Theorem~\ref{mainthmintro} relies on the same idea used in \cite{GWZ20,GWZ201}, which boils down to a Fubini argument along the Hilbert-Chow morphism $h_\chi$ and the fact that up to measure $0$, the fibers of $h_\chi$ over a non-archimedean local field are Picard-schemes of smooth curves.

As a corollary  of Theorem~\ref{mainthmintro}, we obtain from \cite{maulik2020cohomological} (combined with the results in \cite{yuan}) that after pushforward along $h_\chi$, the function $\pBPS_{\beta,\chi}$ agrees with the trace of Frobenius on the $\BPS$-sheaf for the case of del Pezzo surfaces. The $\chi$-independence conjecture for $\BPS$-cohomology suggests, that this continues to hold in the K3 case.

A more concrete application of Theorem~\ref{mainthmintro}, via the Weil conjectures, is the $\chi$-independence of the Betti numbers for $\M_{\beta,\chi}$ as long as $\chi$ and $\beta$ are coprime. 
For $\mathbb P^2$ this is also a consequence of a more general result of Bousseau  \cite[Conjecture~0.4.2,Theorem~0.4.5]{bousseautak}.

When $\M_{\beta,\chi}$ is smooth, the first part of Theorem~\ref{mainthmintro} implies 
\[\pBPS(x)=\BPS(x)\]
for all $x\in\M_{\beta,\chi}(k)$.

 While for more general pairs $(\beta,\chi)$, Definition \eqref{pbpsdef} has no obvious connection with the definition of the $\BPS$ function as trace of Frobenius of some suitable cohomology,  Theorem~\ref{mainthmintro} suggests that the above equality still holds.

For the case of $S$ del Pezzo, this intriguing identity is the subject of an ongoing project of the third named author with Michael Groechenig and Paul Ziegler.

The refined $\chi$-independence for $\BPS$-invariants of moduli of sheaves on a K3 surface is instead still conjectural (see \cite{maulikthomas2019} for a proof of the numerical version), so is the identification of $\pBPS$- and $\BPS$-invariants.

We expect that a conceptual proof (not relying on the $\chi$-independence) for the equality between  $\pBPS$- and $\BPS$-invariants will be easier to obtain in the Fano case, due to the smoothness of the moduli stacks.


Finally, as in \cite{maulik2020cohomological}, we can prove the same statements for moduli of (meromorphic or not)  Higgs bundles; that is pairs of a vector bundle $\mathcal E$ on a curve $C$ of genus $g(C) \geq 2$ and a morphism $\Theta\colon\mathcal E\to\mathcal E(D)$ for $D$ a fixed effective divisor of degree $\deg(D)\geq 2g(C)-2$.

In this case, the Hilbert-Chow morphism is replaced by the Hitchin fibration
  \[H\colon \M_{r,\chi}(C)\to\bigoplus_{i=1}^r \operatorname{H}^0(C,\mathcal O_C(iD))\]
associating to $(\mathcal E,\Theta)$ the chracteristic polynomial of $\Theta$.
The $\chi$-independence for meromorphic Higgs bundles was also proven combinatorially in \cite{MR19}. 
Moreover, during the writing of this article, $\chi$-independence for $\BPS$-invariants of Higgs bundles has been proved in \cite{KK21}. This indirectly implies that our non-archimedean $\pBPS$-function agrees with the geometric $\BPS$-function beyond the Fano cases. 

This is in sharp contrast with the intersection complex on the coarse moduli space, which does depend on the Euler characteristic $\chi$ in the $K3$ and Higgs bundle case, see \cite[Section 0.4]{maulik2020cohomological}.





\subsection*{Acknowledgements} We warmly thank Michael Groechenig and Paul Ziegler for numerous discussions, which were at the origin of many ideas in this paper. We also would like to thank Ben Davison, Pierrick Bousseau, Sven Meinhardt and Tanguy Vernet for interesting conversations around $\BPS$-invariants. We also thank the anonymous referee, whose comments greatly improved the exposition. This work was supported by the Swiss National Science Foundation [No. 196960].

\section{Moduli spaces of sheaves and Higgs bundles}\label{sec:sheaves}

In this section we recall some properties of moduli spaces of sheaves on smooth surfaces and moduli spaces of (meromorphic or not) Higgs bundles relative to a base scheme $T.$

The case $T=\Spec k$ with $k$ algebraically closed of $\operatorname{char}(k)=0$ is classical and all the statements can be found  for example in
\cite[\S~4]{MR1450870} for moduli of sheaves on surfaces, and in  \cite{MR1085642} for moduli of Higgs bundles.

For $T$ a field of positive or mixed characteristic, moduli spaces of sheaves were first constructed by Langer in \cite{langerAnnals,langerDuke}. More general moduli stacks and moduli spaces of sheaves and complexes over a smooth projective family $X/T$ were studied in \cite{lieb} and \cite{huythom}.
Recently the study of moduli stacks in the relative context has been vastly generalised in \cite{stabilityinFamilies}.

 We will often cite the latter reference, even though the results of Langer, Lieblich  and of Huybrechts-Thomas  would suffice for our purposes.
\subsection{Relative moduli spaces of sheaves on surfaces}
Let $S\xrightarrow{r} T$ be a smooth family of projective surfaces, i.e. $r$ is smooth and projective. We will be mostly interested in the case where the relative anti-canonical bundle $-K_{S/T}$ is $r$--ample or trivial; we will then say that $S$ is del Pezzo, respectively K3, over $T.$ 

 For us $T$ will be either $\Spec (\mathcal O)$, or
$\Spec F,$ or $ \Spec  k$ where $F$ is a p-adic field, $\mathcal O$ its ring of integers, and $k$ a field, possibly of positive characteristic.
\subsubsection*{Construction of the moduli spaces}
 Let us fix $\beta$ on $S$ the class of a $r$--effective divisor, ample and base point free, i.e. $\beta_{t}$ is ample and base point free for each $t \in T$. Given $\mathcal F$ a sheaf on $S$, flat over $T$, we say that $\mathcal F$ is \emph{1-dimensional} over $T$ of class $\beta$ if \[\operatorname{Supp}(\mathcal F):=\operatorname{Fit}_0(\mathcal F)= C\]
 with $C\in \mathbb P_T(g_*\mathcal{O}_S(\beta)).$ Since the Fitting support is compatible with pull-back \cite[\href{https://stacks.math.columbia.edu/tag/07Z6}{Tag 07Z6}]{stacks-project} this means that $ \mathcal F_{t}$ is supported on a curve of class $\beta_{t}$ on $S_{t}$ for each $t \in T$.

Let $\mathcal O_S(1)$ be a $r$-polarization. Given $\mathcal F$ a sheaf on $S$ flat over $T,$  the Hilbert polynomial
\[P_r(n):=\chi(S_{t}, \mathcal F_{t}\otimes\mathcal O_{S_t}(n))=(\beta_{t}\cdot H_{t})n +\chi(\mathcal F_{t})\]
 for $ H_{t}\in |\mathcal O_{S_t}(1)|$
is locally constant on $T$ (see for example \cite[Theorem~9.9]{MR0463157}); in particular $\chi(\mathcal F_{t})$ is locally constant.

 We denote by $p_r(n)$ the reduced \emph{relative Hilbert polynomial} obtained from $P_r(n)$ dividing by the coefficient of the leading term.

\begin{definition}(see  \cite[Definition~21.11]{stabilityinFamilies} )

Let $\beta,\chi$ and $\mathcal O_S(1)$ be as above. We denote by
\[\mathbb M_{\beta,\chi}(S)\colon(\operatorname{Sch}/T)^{\text{op}}\rightarrow\operatorname{Gpds}\]
the functor whose value on $W\to T$ consist of sheaves $\mathcal F$ on $S\times_T W,$ flat over $W$ such that for all $w\in W$, $\mathcal F_w$ is a Gieseker semi-stable  1-dimensional sheaf of class $\beta_w$ and Euler characteristic $\chi\in\mathbb Z.$
We recall that $\mathcal{F}_w$ is said Gieseker semi-stable if for all proper subsheaves $\mathcal G\in\mathcal{F}_w$  we have:
\begin{equation}\label{eq:stability}
n +\frac{\chi(\mathcal G)}{[\operatorname{Supp}(\mathcal F_w)]\cdot H_{w}} \leq p_r(n),
\end{equation} 
and Gieseker stable if the inequality is strict. Finally we say that $\mathcal{F}_w$ is geometrically Gieseker stable if its pull-back  $\mathcal{F}_{\overline{w}}$ to the algebraic closure is Gieseker stable.

An object $\mathcal F\in \mathbb M_{\beta,\chi}(S)(W)$ is called a \emph{family of  Gieseker semi-stable} sheaves of class $(\beta,\chi)$ over $W.$

$\mathbb M^{\text{st}}_{\beta,\chi}(S)$ is the sub-functor of families of Gieseker geometrically stable objects.

\end{definition}

\begin{theorem}\label{thm:semistablesheaves}\cite{langerAnnals,langerDuke,stabilityinFamilies}
In the notation above:
\begin{enumerate}
\item $\mathbb M_{\beta,\chi}(S)$ and $\mathbb M^{\text{st}}_{\beta,\chi}(S)$ are algebraic stacks of finite type over $T.$ The same statement hold for the $\mathbb G_m$-rigidification $\mathcal M_{\beta,\chi}(S)$ of $\mathbb M_{\beta,\chi}(S).$ 
\item There exists a projective $T$-scheme $\M_{\beta,\chi}(S)$ which uniformly represent the set valued functor over $T$-schemes associating to $W$ the set of $\text{S}$-equivalence classes of families of Gieseker semi-stable sheaves on the fibers of $S\times_T W\to W$ with Hilbert polynomial $p_r(n).$
\item $\mathbb M^{\text{st}}_{\beta,\chi}(S)$ is a $\mathbb{G}_m$-gerbe over its coarse moduli space $\M^{\text{st}}_{\beta,\chi}(S)$, which is a quasi-projective scheme over $T$ representing isomorphism classes of  Gieseker geometrically stable sheaves on the fibers of $S\times_T W\to W$ with Hilbert polynomial $p_r(n).$
\item  $\mathbb M_{\beta,\chi}(S)$ is universally closed over $T$ and thus universally closed over $\M_{\beta,\chi}(S).$ 

\end{enumerate}

\end{theorem}
\begin{proof}

(1) A far more general version of this statement is proved in \cite[Theorem~21.24]{stabilityinFamilies}.
  For the benefit of the reader, we recall the key steps of the proof in the special case of Langer \cite{langerDuke}, which is the only one we need. 
  
Let $E$ be a coherent sheaf on $S$ flat on $T$, and let  $\operatorname{Quot}_{T,p_r(n)}(E)\colon (\operatorname{Sch}/T)^{\text{op}}\rightarrow\operatorname{Sets}$ be the relative Quot functor, whose $W\to T$ valued points are quotients $E_W\to Q\to 0$ with $Q$ a coherent sheaf on $S\times_T W$, flat over $W$, with Hilbert polynomial $p_r(n)$ along the fibers of  $S\times_T W\xrightarrow{r_W} W$. As proved for example in \cite[\href{https://stacks.math.columbia.edu/tag/09TQ}{Tag 09TQ}]{stacks-project}, $\operatorname{Quot}_{T,p_r(n)}(E)$  is represented by an algebraic space over $T$ locally of finite presentation. 

By Langer's results \cite[\S~3]{langerDuke}, families of Gieseker semi-stable sheaves on the fibers of $S\xrightarrow{r} T$ with fixed Hilbert polynomial $p_r(n)$ are bounded, i.e. any $\mathcal F\in \mathbb M_{\beta,\chi}(S)$ is a quotient of $E=\mathcal O_S(-N)^{\oplus F(N)}$ for some fixed $N.$

Furthermore, as showed in \cite[\S~4]{langerDuke}, the open subsets $\mathcal R^{\text{st}}\subseteq \mathcal R^{\text{ss}}\subset \operatorname{Quot}_{T,p_r(n)}(E)$  of Gieseker stable/ semi-stable sheaves on the fibers of $S\xrightarrow{r} T$, correspond to the opens of stable/ semi-stable points with respect to the natural action of the $T$-smooth group $\operatorname{GL}(F(N)).$

It thus follows that 
$\mathbb M_{\beta,\chi}(S)\cong [ \mathcal R^{\text{ss}}\slash \operatorname{GL}(F(N))]$ and $\mathbb M^{\text{st}}_{\beta,\chi}(S)\cong[ \mathcal R^{\text{st}}\slash \operatorname{GL}(F(N))]$ as stacks over $\operatorname{Sch}/T.$
The statement for the rigidification $\mathcal M_{\beta,\chi}(S)$ now follows immediately as the latter is simply given by $[ \mathcal R^{\text{ss}}\slash \operatorname{PGL}(F(N))].$

(2) The projective $T$-scheme $\M_{\beta,\chi}(S)$  is the GIT quotient $\mathcal R^{\text{ss}}\slash\slash \operatorname{GL}(F(N))$ this exists and has the claimed property by \cite[Theorem~4]{Seshadri}, as also argued by Langer in the proof of \cite[Theorem~4.1]{langerDuke}.

(3) The statement follows by the construction explained in $(1)$.

(4) The universal closedness over $T$ follows from a relative version of Langton's Theorem. The reader can see \cite{MR0364255} or \cite[Theorem~2.B.1]{MR1450870} for the proof over an algebraic closed field, and \cite[Lemma~21.22]{stabilityinFamilies} for the general case.
 
Given the GIT construction of the moduli space recalled above, the universal closedness of  $\mathbb M_{\beta,\chi}(S)\to \M_{\beta,\chi}(S)$ also follows from \cite[Theorem~A.8]{AHH18}.

\end{proof}
\subsubsection*{Deformation Theory}
Let $S\xrightarrow{r} T$ be as in the previous section, and let $\mathcal F$ be a coherent sheaf on $S$ flat over $T.$ 
We denote by  $\operatorname{Ext}^i_S(\mathcal F,\mathcal F)$ the $T$-sheaf obtained by taking the i-th cohomology of the complex $Rr_*\R\mathcal Hom(\mathcal F,\mathcal F),$ i.e.:
\[\operatorname{Ext}^i_S(\mathcal F,\mathcal F):= \text{H}^i(Rr_*\R\mathcal Hom(\mathcal F,\mathcal F)).\]

By \cite[\S~3]{lieb} (see also \cite[\S~4]{huythom}) the sheaves $\operatorname{Ext}^i_S(\mathcal F,\mathcal F)$
control the deformation theory of $\mathcal{F}$ over $T$, namely:
\begin{enumerate}
\item There is a natural isomorphism  $\operatorname{Ext}^0_S(\mathcal F,\mathcal F)\cong \operatorname{End}(\mathcal F,\mathcal F).$
\item The first order deformations of $\mathcal F$ are given by  $\operatorname{Ext}^1_S(\mathcal F,\mathcal F)$
\item The obstruction to deformations are contained in $\operatorname{Ext}^2_S(\mathcal F,\mathcal F).$
\end{enumerate}

\begin{lemma}\label{lem:relativesmooth}
Let $S\xrightarrow{r} T$ be a smooth family of del Pezzo surfaces. Then $\BM_{\beta,\chi}(S)$ is smooth over $T.$
\end{lemma}
\begin{proof}
It suffices to show that the obstruction $\operatorname{Ext}^2_S(\mathcal F,\mathcal F)$ is zero. We will argue that $\operatorname{Ext}^2_S(\mathcal F,\mathcal F)\otimes k(\bar{t})=0$ for any geometric point $\bar{t}\in T.$
By \cite[Theorem~1.9]{altman}, the cohomology and base change property for Ext sheaves implies:
\[\operatorname{Ext}^2_S(\mathcal F,\mathcal F)\otimes k(\bar{t})=\operatorname{Ext}^2_{S_{\bar{t}}}(\mathcal F_{\bar{t}},\mathcal F_{\bar{t}}),\]
where $\mathcal F_{\bar{t}}$ is by definition a Gieseker semi-stable sheaf of class $\beta_{\bar{t}}$ and Euler characteristic $\chi$ on $S_{\bar{t}}.$
The vanishing of the RHS is then classical, see for example \cite[Lemma~2.5]{maulik2020cohomological}.
\end{proof}
\subsubsection*{Singularities in the K3 case}
For $S\to T$ a K3 surface, it will never be the case that $\operatorname{Ext}^2_S(\mathcal F,\mathcal F)$ vanishes. Nonetheless, the following generalization of a classical result of Mukai \cite{mukai1984} holds:

\begin{theorem}\label{thm:splK3}\cite[Theorem 3.2]{inaba2011smoothness}
Let $S\slash k$ be a smooth projective K3 surface over an algebraically closed field $k$ and let $\mathcal F$ be a simple sheaf, i.e. $\operatorname{Hom}_{S}(\mathcal F,\mathcal F)=k.$ Then the deformations of $\mathcal F$ are unobstructed. In particular, the open sub-stack $\operatorname{Spl}(S\slash k)_{(\beta,\chi)}\subseteq \mathbb M_{\beta,\chi}(S)$ of simple semi-stable sheaves is smooth over $k.$
\end{theorem}

\begin{corollary}\label{k3cor} Let $S \to T$ be a smooth projective relative K3 surface and $T$ a smooth base with perfect residue fields. Then the moduli stack $\operatorname{Spl}_{(\beta,\chi)}(S\slash T)$ of simple semi-stable sheaves is smooth over $T$.
\end{corollary}
\begin{proof}
 By Inaba's result, Theorem~\ref{thm:splK3}, the closed fibers of $\operatorname{Spl}_{(\beta,\chi)}(S\slash T)\to T$ are smooth of dimension $\beta_{\bar{t}}^2+2,$ and since $S\to T$ is smooth this number is constant on $T$ \cite[\S~20.3]{fultonintersection}.
  It then follows from \cite[Theorem~3.3.27]{schoutens2010use} that $\operatorname{Spl}_{(\beta,\chi)}(S\slash T)\to T$ is flat and thus in fact smooth since all of the closed fiber are.
\end{proof}

This immediately implies that $\mathbb M^{\text{st}}_{\beta,\chi}(S)$ is smooth for $S \to T$ as in Corollary \ref{k3cor}. Furthermore, for $T=\mathbb C$, the singularities of the moduli stack $\mathbb M_{\beta,\chi}(S)$ have been intensively studied \cite{arbarellosacca,budur2019formality,davison2021purity,crawley2001geometry}.
In the following theorem, we state, in a form convenient for our purposes, some consequences of the results obtained in the quoted references.

\begin{theorem}\label{thm:K3case}
Let $S\slash\mathbb C$ be a smooth projective K3 surface and $(\beta,\chi)$ a dimension vector such that the locus of $\operatorname{Spl}_{\beta,\chi}^{\text{ss}}(S)\subseteq\mathbb M_{\beta,\chi}(S)$ of simple semi-stable sheaves has complement of codimension at least two. Then $\mathbb M_{\beta,\chi}(S)$ is normal with local complete intersection singularities.
\end{theorem}
\begin{proof}

For a morphism $Y\to X$ he property of being local complete intersection is \'etale local on source and target  \cite[\href{https://stacks.math.columbia.edu/tag/06C3}{Tag 06C3}]{stacks-project}. Thus it is enough to argue that for any point $q\in  \mathbb M_{\beta,\chi}(S)$ there exist an \'etale neighbourhood $U_q$ with l.c.i. singularities. 

By \cite[Theorem~5.11]{davison2021purity}, an \'etale local model for $\mathbb M_{\beta,\chi}(S)\to \text{M}_{\beta,\chi}(S)$ is given by the moduli stack $\mathfrak M_{\mathbf{d}}(\Pi_Q)=[\mu_{\mathbf d}^{-1}(0)/\GL_{\mathbf{d}}]$  (see \cite[\S~3]{davison2021purity} for notation) of representations of the preprojective algebra of a certain explicit quiver, already considered in \cite{arbarellosacca}.

On the other hand, by \cite[Theorem~1.2]{crawley2001geometry}, $\mu_{\mathbf{d}}^{-1}(0)$ is a local complete intersection as soon as there exists a simple representation of the preprojective algebra with given dimension vector. 

By our assumptions, $\operatorname{Spl}_{\beta,\chi}^{\text{ss}}(S)$ is dense in $\mathbb M_{\beta,\chi}(S)$. This ensures the existence of a simple representation in $U_q$ for each $q$. Thus, since étale morphisms preserve the dimension of the automorphism groups and automorphism groups for quiver representations are connected, there exists a simple representation for $\Pi_Q$ for each $Q$ appearing in the local description.

Once we know that $\mathbb M_{\beta,\chi}(S)$  is l.c.i, and thus Gorenstein, the normality simply follows from the assumption on the codimension of the locus of simple sheaves and Mukai-Inaba's Theorem \ref{thm:splK3}. \end{proof}



\begin{rmk}\label{coas}
We show in Section~\ref{sec:invertiblesheaves} that the hypothesis on the codimension of the complement of $\operatorname{Spl}_{\beta,\chi}^{\text{ss}}(S)$ is satisfied as soon as: $\beta$ is ample and base point free, and the locus $V$ of non-reduced curves in the linear system $|\beta|$ is of codimension  at least 2.

Let for example $S$ be the K3 surface obtained as the 2:1 cover of $S\xrightarrow{\tau}\mathbb P^2$ branched along a general sextic $Z.$ Consider $\beta=\tau^*[dL]$ for $d\geq 1.$ Then a curve $C$ in  $|\beta|$ is the 2:1 cover of a curve $D\in |dL|$ branched over $D\cap Z.$ Then $C$ can be non-reduced only if, either $D$ itself is non-reduced or there is a component of $D$ contained in $Z,$ i.e. $Z$ is a component of $D.$ Both loci have high codimension in  $|dL|.$
\end{rmk}

\subsubsection*{Hilbert-Chow map}
We keep the notation of the previous section. 
Let us denote by $B_T=\mathbb P_T(g_*\mathcal O_S(\beta)),$ and by $\mathcal  C\to B_T$ the universal curve.
As for the case $T=\Spec(\mathbb{C})$ there is a morphism
\[\mathbf{h}_{\chi}\colon \mathbb M_{\beta,\chi}(S)\to B_T,\]
which associate to $\mathcal F$ its fitting support and factors through the moduli space
\[h_{\chi}\colon \M_{\beta,\chi}(S)\to B_T.\]
We call $h_{\chi}$ the Hilbert-Chow morphism and $\mathbf{h}_{\chi}$ the stacky Hilbert-Chow morphism. 
Since both $ \M_{\beta,\chi}(S)$ and $B_T$ are projective over $T,$ $h_{\chi}$ is proper over $T.$

For $b\colon T\to B_T$ a section such that $\mathcal C_b/T$ is smooth, the fiber $\mathbf{h}_{\chi}^{-1}(b)$ is the Picard stack $\mathfrak{P}ic^{\chi+g-1}_{\mathcal C_b/T}$ where $g$ is the genus of the curve, and $h^{-1}_{\chi}(b)$ the Picard space $\operatorname{Pic}^{\chi+g-1}_{\mathcal C_b/T},$ i.e. the $\mathbb G_m$-rigidification.

For $b\colon T\to B_T$, such that $\mathcal C_b/T$ is integral, the fiber of the stacky Hilbert-Chow morphism is some universally closed stack over $T$ containing  $\mathfrak{P}ic^{\chi+g-1}_{\mathcal C_b/T}$ as an open sub-stack, and $h^{-1}_{\chi}(b)$ is the corresponding GIT quotient, which is the compactification of  the Picard space $\operatorname{Pic}^{\chi+g-1}_{\mathcal C_b/T}$ constructed by Altman-Kleiman \cite{altman}.

In general, for $\mathcal C_b/T$ non-integral not all line bundles of total degree $\operatorname{totdeg}=\chi+g-1$ are Gieseker semi-stable and one cannot say much about the fiber, except that it is a moduli stack of semi-stable sheaves of pure dimension 1.
In the reduced case, moduli stacks and moduli spaces of rank 1 semi-stable torsion free sheaves on  $\mathcal C_b/T$ have been studied, and below we will recall some known results.
\subsection{The open sub-stack of invertible sheaves}\label{sec:invertiblesheaves}
We denote by
\[\mathfrak{U}_{\beta,\chi}\subseteq\mathfrak Pic^{\chi+g-1}_{\mathcal C/B_T}\]
the open sub-stack
 corresponding to line bundles of total degree $\chi+g-1$, which are Gieseker semi-stable on the fibers of $\mathcal C/B_T$.
 
The rigidification $ \mathcal U_{\beta,\chi}\subseteq \mathcal M_{\beta,\chi}(S)$ is stabiliser free, and 
$\mathfrak{U}_{\beta,\chi}\to \mathcal U_{\beta,\chi}$ is simply a $\mathbb G_m$--gerbe. It follows from Theorem \ref{thm:semistablesheaves}, that if we further restrict to geometrically stable line bundles $\mathfrak{U}^{\text{st}}_{\beta,\chi}$, the rigidification  $\mathcal U^{\text{st}}_{\beta,\chi}$ is a good moduli space for the moduli stack and it is isomorphic to an open in $\M_{\beta,\chi}(S).$

The opens $\mathcal U^{\text{st}}_{\beta,\chi}\subseteq\mathcal U_{\beta,\chi}$ will play a crucial role in the construction of the canonical measure.
To this end, we need to prove the following codimension estimate:
\begin{proposition}\label{prop:codim}
\begin{enumerate}
\item Let $S\xrightarrow{r} T$ be a relative del Pezzo surface and $\beta$ a base point free and ample curve class.
Then the complement of $\mathcal {U}_{\beta,\chi}$ in $ \mathcal M_{\beta,\chi}(S)$, or equivalently the complement of  $\mathfrak {U}_{\beta,\chi}$ in $\BM_{\beta,\chi}(S)$, has codimension at least two.
\item Let $T=\operatorname{Spec}(\mathbb C)$ and $S$ a smooth projective K3 surface. Suppose furthermore that $\beta$ is ample, base point free and the codimension of the non-reduced locus $V\subset B$ is at least 2. Then the complement of $\mathcal {U}_{\beta,\chi}$ in $ \mathcal M_{\beta,\chi}(S)$ has codimension at least two.
\end{enumerate}

\end{proposition}
The proof of this statement is somewhat involved and it will occupy the rest of the section. The reader who wish to get to the core of the paper quickly may skip this in a first reading.

\begin{proof}[Proof of~\ref{prop:codim}]
In case $(1)$, by Lemma~\ref{lem:relativesmooth}, $\BM_{\beta,\chi}(S)\to T$ is smooth, so it is sufficient to estimate the codimension of the locus of line bundles on each geometric fiber $\BM_{\beta_{\bar{t}},\chi}(S_{\bar{t}})\to \Spec(\bar{k})$. 
In case $(2)$ we are already by hypothesis working on an algebraically closed field.
This allows us to reduce to the following classical setting:
\begin{setting}
$S$ is a smooth del Pezzo or K3 surface over $\bar{k}$ an algebraically closed field; $H$ is a polarization on $S$, $\beta$ is a base point free and ample curve class, and $B=|\beta| $ the linear system of curves in class $\beta.$
\end{setting}

Working on $T=\Spec(\bar{k})$ allows us to appeal to the results of \cite{maulik2020cohomological,yuan}. It is proved in \cite{yuan}, with techniques inspired by \cite{schiffmann}, that for smooth projective surfaces $S$ over an algebraically closed field $\bar{k}$ with canonical $K_S\leq 0$ the stacky Hilbert-Chow morphism $\mathbf{h}_{\chi}$ (and thus also $h_{\chi}$) is equidimensional:
\begin{theorem}\cite[Theorem~1.2, Corollary~1.3]{yuan}\label{Thm:equidimensionality}
Let $S$ be a smooth projective surface and $\beta=\sum_{i=1}^n m_i\beta_i$ an ample effective curve class with $\beta_i$ pairwise distinct. If $K_S\cdot\beta_i\leq 0$ for all $i$, then 
\[\dim(\mathbf{h}^{-1}_{\chi}(b))\leq g-1\]
for any $b\in |\beta|$, where $g=\frac{1}{2} \beta\cdot(\beta+K_S)+1.$ 
\end{theorem}
Using the dimension estimate, Maulik-Shen prove the following:

\begin{corollary} \cite[Theorem~2.3]{maulik2020cohomological}\label{cor:irr}
For $S$ del Pezzo, Hilbert-Chow morphisms $\mathbf{h}_{\chi}$ and $h_{\chi}$ are equidimensional.
In particular $\M_{\beta,\chi}(S)$ is irreducible. 
\end{corollary}
\begin{corollary}
For $S$ a K3 over an algebraically closed field, Hilbert-Chow morphisms $\mathbf{h}_{\chi}$ and $h_{\chi}$ are equidimensional.
In particular $\M_{\beta,\chi}(S)$ is irreducible. 
\end{corollary}
\begin{proof}
This is a word by word repetition of the argument in  \cite[Theorem~2.3]{maulik2020cohomological} given the dimension estimate in Theorem~\ref{Thm:equidimensionality}.
\end{proof}
\begin{rmk}
The dimension estimate $\dim(\mathbf{h}^{-1}_{\chi}(b))\leq g-1$ had previously been proved in  \cite[Proposition~2.6]{maulik2020cohomological} under the additional hypothesis that $S$ is a toric del Pezzo, by reducing the estimate to curves with underlying support a union of possibly non-reduced toric divisors, for which the results of \cite{Laumonsupp} apply.
\end{rmk}
We use the equidimensionality of the (stacky) Hilbert-Chow as follows. First, by Bertini's Theorem, since $\beta$ is ample and base point free, there is a dense open $B^{sm}\subset B$ such that for $b\in B^{sm}$ the curve $\FC_b$ is a non-singular, integral curve of genus  $g=\frac{1}{2} \beta\cdot(\beta+K_S)+1.$  In this case, $\mathbf{h}_{\chi}^{-1}(b)\cong\mathfrak Pic^{\chi+g-1}_{\FC_b},$ which is simply a $\mathbb G_m$-gerbe over the degree $\chi+g-1$ component of the Picard group.

In particular $\mathbf{h}_{\chi}^{-1}(B^{sm})$ is contained in the open sub-stack of (stable) line bundles. 
Notice furthermore that since the number of connected components is lower semi-continuous in a proper flat family (see for example \cite[Theorem~4.17]{deligne1969irreducibility}) all the curves in the linear system are connected.

To prove the codimension estimate, we now look at the fibers of the Hilbert-Chow morphism over the locus of singular curves.

The complement of $B^{sm}$ in $B$ has codimension 1, and by  equidimensionality of $\mathbf{h}_{\chi}$ so does $\mathbf{h}^{-1}_{\chi}(B\setminus B^{sm}).$ 

Let us denote by $V\subseteq B\setminus B^{sm}$ the locus of non-reduced curves.  We will prove in Section~\ref{sec:nonred} that this has always codimension at least 2 in $B$ for $S$ del Pezzo,  provided  $\beta$ is a base point free ample class. For $S$ a K3, the codimension of $V$ is at least two by assumption; see Remark~\ref{coas} for examples where the assumption is verified. Then $\mathbf{h}^{-1}_{\chi}(V)$  has codimension at least two, and it can be discarded for the purposes of estimating the codimension of the complement of $\mathfrak{U}_{\beta,\chi}.$

Finally, we study the fibers over $(B\setminus B^{sm})\setminus V$ in Section~\ref{sec:reducedsing} and prove that the locus of line bundles in $\mathbf{h}_{\chi}^{-1}(b)$ is dense for all $b\in (B\setminus B^{sm})\setminus V.$ 

In conclusion, the complement of  $\mathfrak{U}_{\beta,\chi}$ is contained in the union of $\mathbf{h}^{-1}_{\chi}(V)$ and a divisor $\mathfrak{D}\subset \mathbf{h}^{-1}_{\chi}(B\setminus B^{sm}).$ These have both codimension at least two, which allows us to conclude.

\subsubsection{Excluding the non-reduced locus}\label{sec:nonred}
 Consider the stratification of $B$ given by
\[B=\bigsqcup_{\underline{\beta}} B_{\underline{\beta}},\]
where:
 $\underline{\beta}=\left\{m_1\beta_1,\dots ,m_s\beta_s\right\}$ with $s\geq 1$, $m_i\geq 1$, and $\beta_i$ are distinct curve classes such that there exists an integral curve in $|\beta_i|$ for each $i$;
 $B_{\underline{\beta}}$ is the locally closed strata in $B$ of curves of type $\underline{\beta}$, i.e.
\[B_{\underline{\beta}}=\left\{ C=\sum_i m_i C_i\;\in\;|\beta|,\; \; C_i\;\in\;|\beta_i|\right\},\]
with $C_i$ integral.

Let $Z$  be an irreducible component in $\overline {B}_{\underline{\beta}}$; we denote by $\delta_Z$ the dimesion of the affine part of the Jacobian of $\mathcal C_b$ for  $b\in Z^{\circ}\subseteq \overline {B}_{\underline{\beta}}$ generic, i.e.
\[\delta_Z=\dim(\operatorname{Pic}^0(\mathcal C_b))^{\text{Aff}}.\]

\begin{lemma}\label{lem:codimnonred}
Let $\beta$ be a base point free ample class on a smooth del Pezzo surface. The locus $V\subset B$ of non-reduced curves has codimension at least two.
\end{lemma}
\begin{proof}
The locus $V$ of non-reduced curves is the union of closed strata $\overline {B}_{\underline{\beta}}$ with $\underline{\beta}$ having at least one $m_i\geq 2.$ 
We want to show that the codimension of  any irreducible component $Z$ of such a stratum $\overline {B}_{\underline{\beta}}$ is at least two in $B$.
The proof will follow from the following inequality which is proved in \cite[Proposition~4.3]{maulik2020cohomological}:
\begin{equation}\label{eq:MSinequality}
\operatorname{codim}(Z)\geq -K_S\cdot (\sum_i(m_i-1)\beta_i)+\delta_Z\geq 1+\delta_Z
\end{equation}
where the second inequality holds since  $-K_S$ is ample, $m_i-1\geq 0$ for each $i$ and $m_i-1\geq 1$ for at least one $i$ and thus
\[-K_S\cdot (\sum_i(m_i-1)\beta_i)\geq 1;\]
since $-K_S$ is ample, the equality can hold only if all but one $m_i$ are equal to $1$.

If the latter inequality is  strict, we are done; else we can assume that $m_1=2$ and $m_i=1$ for any other $i$.
Again,  if $\delta_Z\geq 1$ we are done.

To complete the proof we only need to look the cases where  $m_1=2$, $m_i=1$ for any other $i$, and  $\delta_Z=0$. 

Let us denote by $C:=\mathcal C_b$ for $b\in B_{\left\{2\beta_1,\beta_2,\dots,\beta_s\right\}}$ a generic point in this stratum, by $C'\subseteq C$ the reduced subcurve, and by $\beta'=\beta-\beta_1$ its class.

We have  two short exact sequences: 
\begin{align*}
0\to \mathcal O_{C_1}(-\beta'\cdot\beta_1)\to \mathcal O_{C}\to\mathcal O_{C'}\to 0\\
 0\to O_{C'}\to \oplus\bigoplus_{i=1}^s\mathcal O_{\widetilde{C}_i}\to \bigoplus_{p\in C'^{\text{sing}} } \mathcal O_p^{\delta(p)}\to 0.
  \end{align*}
The latter is simply the normalisation sequence for $\widetilde{C'}=\sqcup_{i=1}^s\widetilde{C}_i\to C'$ (see \cite[\S~7.5]{Liu}). We denoted by 
\[\delta(p)=\operatorname{length}_k(\mathcal O_{\widetilde{C'},p}\slash\mathcal O_{C'})\]
the so called delta invariant of the singularity.

By \cite[Lemma~5.11]{Liu}, there is a surjective morphism of abelian groups
\[\operatorname{Pic}^0(C)\to \operatorname{Pic}^0(C')\to 0\]
with kernel an affine unipotent group of dimension $\operatorname{dim}\operatorname{H}^1(C_1,\mathcal O_{C_1}(-\beta'\cdot\beta_1)).$
Moreover, by \cite[Theorem~5.19]{Liu}, $\operatorname{Pic}^0(C)$ is an extension of $\times \operatorname{Pic}^0(\widetilde{C}_i)$ by an affine group of dimension 
\begin{equation}\label{eq:deltared}
\left( \sum_{p\in C'^{\text{sing}}} \delta (p)\right) -s+1,
\end{equation}
 and thus
\[\delta_Z= \left( \sum_{p\in C'^{\text{sing}}}{ \delta (p)}\right) -s+1+\operatorname{dim}\operatorname{H}^1(C_1,\mathcal O_{C_1}(-\beta'\cdot\beta_1)).\]

One can prove by induction on the number $s$ of irreducible components that $\sum_{p\in C'^{\text{sing}}} \delta (p) -s+1\geq 0$ with equality holding only if the only singularities are nodes and $C'$ is of compact type, namely its dual graph is a tree.
For the induction step look at the partial normalisation of $C'$ given by $R\sqcup\overline{C'\setminus R}$ for $R$ an irreducible component.

Let us also analyse the degree of $\mathcal O_{C_1}(-\beta'\cdot\beta_1).$ We can write $\beta'=\beta"+\beta_1$ with $\beta"$ and effective curve class, and we can assume that $\beta"\cdot\beta_1>0$ otherwise $C'$ is not connected and we can estimate $\delta_Z$ separately on each connected component.
From $0\leq g(\beta_1)=\frac{1}{2} K_S\cdot\beta_1+  \frac{1}{2} \beta^2_1+1$ and the fact that $K_S$ is anti-ample, it follows that $\beta^2_1\geq -1.$ In particular, 
$-\beta'\cdot\beta_1=- \beta"\cdot\beta_1-\beta^2_1\leq 0$ with equality holding only in the case $\beta_1^2=-1$, and $\beta"\cdot\beta_1=1.$ However, if the equality is verified, then $\beta\cdot\beta_1 <0$, contradicting the ampleness of $\beta.$

We can thus assume that  $\mathcal O_{C_1}(-\beta'\cdot\beta_1)$ has strictly negative degree $d.$ By Riemann-Roch, 
\[\operatorname{dim}\operatorname{H}^1(C_1,\mathcal O_{C_1}(-\beta'\cdot\beta_1))=g(C_1)-d-1,\]
which is always positive unless $g(C_1)=0$ and $d=-1.$ In particular in this case we must have:  $\beta"\cdot\beta_1=1,\beta_1^2=0,$ and $K_S\cdot\beta_1=-2$
We can thus give a better estimate of the codimension of the component $Z$ in $B$ corresponding to such curves. By \cite[Lemma~2.1]{maulik2020cohomological}, given $\gamma$ an effective curve class, the dimension of the linear system $|\gamma|$  is $\frac{1}{2}\gamma\cdot(\gamma-K_S).$
 So $\dim B=\frac{1}{2}\beta\cdot(\beta-K_S),$ while the dimension of $Z$ is $\frac{1}{2}\beta'\cdot(\beta'-K_S).$ We compute:
 \begin{align*}
 \operatorname{codim}(Z)&=\frac{1}{2}\beta\cdot(\beta-K_S)-\frac{1}{2}\beta'\cdot(\beta-K_S)\\
 &=\frac{1}{2}[(\beta'+\beta_1)^2 -K_S\cdot\beta-\beta'^2+K_S\cdot\beta')\\
  &=\frac{1}{2}[2\beta'\cdot\beta_1-K_S\cdot\beta_1]=2
 \end{align*}

\end{proof}

The assumption in case $(2)$ and Lemma~\ref{lem:codimnonred} in case $(1)$ of Proposition~\ref{prop:codim} allows us to disregard completely the locus $V$ of non-reduced curve.  

\subsubsection{Reduced curves}\label{sec:reducedsing}

 If $\FC_b$ is integral, then any rank 1 torsion free sheaf $\mathcal F$ is strictly stable, independently from the choice of polarization. This is easily seen noticing that for any proper sub-sheaf $\mathcal{G}\hookrightarrow \mathcal F$, the cokernel is supported on a dimension $0$ sub-scheme, which imply the strict inequality in \eqref{eq:stability}.
In particular,   $\mathbf{h}_{\chi}^{-1}(b)$ is a $\mathbb G_m$-gerbe over the compactification $\overline{\Pic}^{\chi+g-1}_{\FC_b}$of $\Pic^{\chi+g-1}_{\FC_b}$ of Altman-Kleiman \cite{altman}. It was proved in \cite[Theorem~A]{MR584085} that for an integral curve with only planar singularities $\overline{\Pic}^{\chi+g-1}_{\FC_b}$ is irreducible and thus the open of line bundles is dense.

We are left to study the fiber $\mathbf{h}^{-1}_\chi(b)$ for $b \in |\beta|$ a reduced but possibly reducible curve. 

\begin{proposition}\label{prop:linebundlestack}
Let $C$ be a reduced curve with planar singularities and polarisation $H_C.$ Let $\BM_{\chi}(C)$ be the moduli stack of Gieseker semi-stable rank one torsion free sheaves on $C$. Then the open substack $\mathfrak{U}_{\chi}$ of semi-stable invertible sheaves is dense in $\BM_{\chi}(C)$.
\end{proposition}

\begin{proof}

 To show that $\mathfrak{U}_{\chi}\subseteq \BM_{\chi}(C)$ is dense, we need to show that any point $F$ in $\BM_{\chi}(C)$ is the limit of a family of line bundles on $C$ in a one parameter family. This now follows from the deformation theory of rank $1$ sheaves on planar curves.
Let  $q_1,\dots q_k$ be the points (necessarily) in the singular locus $ C^{\text{sing}}$ where $F$ fails to be a line bundle. We look at the forgetful morphism
 \[\operatorname{Def}_{F}\xrightarrow{l}\prod_{i=1,\dots k} \operatorname{Def}_{F_{q_i}},\]
 where $F_{q_i}$ is the stalk of $F$ at $q_i$. Since we are looking at curves with planar singularities, the forgetful morphism is smooth \cite{fantechi}. Thus it is sufficient to show that the generic element in each local deformation space $R_{F_{q_i}}$ corresponds to a invertible module. But this follows from the analogous statement for compactified Jacobians of integral planar curves \cite{MR584085}. More explicitly: in order to describe $R_{F_{q_i}}$, consider an integral curve $C_{q_i}$ with $q_i$ as unique singular point. By  \cite{MR584085} the Picard is dense in the compactified Jacobian, which in particular implies that the generic point in $R_{F_{q_i}}$ parametrises a line bundle.
\end{proof}

This concludes the proof of Proposition~\ref{prop:codim}.
\end{proof}

\begin{rmk}
When there are no strictly semi-stables, and thus $\BM_{\chi}(C)$ is a $\mathbb G_m$-gerbe on $h_{\chi}^{-1}(C)=\M_{\chi}(C)$, it is stated in \cite[Corollary~2.20]{melorapagnettaviviani} or also in \cite[Theorem~4.5]{LopezMartin} that stable line bundles are dense in the moduli space.
 
When there are strictly semi-stables (which is the most interesting case for us) the density statement at the \emph{level of moduli spaces} $\M_{\chi}(C)$ is false; the description of the moduli space in \cite[Theorem~4.5]{LopezMartin} explains this fact.

\end{rmk}

\subsection{Moduli of Higgs bundles}
Let $C\to T$ be a smooth projective curve. As for the case of sheaves on surfaces, there is an adaptation of theory of moduli stack and moduli spaces of Higgs bundles relative to a base $T.$ The case where $T$ is a field of positive characteristic is already  considered in \cite{MR2218781}.

We refer to \cite[\S~7.1]{GWZ20} and references therein for the adaptation of the classical theory to the case of smooth curves over $T=\Spec (O_F),\Spec(F).$

We denote by $\BM_{r,\chi}(C)$ and  $\M_{r,\chi}(C)$ the moduli stack and moduli space of semi-stable Higgs bundles $(\mathcal E, \Theta\colon\mathcal E\to\mathcal E\otimes\mathcal O_C(D))$ of rank $r$ and Euler characteristic $\chi$ on the fibers of $C\to T.$
Here $D$ is a relative effective divisor of degree $d=\text{deg}(D) \geq 2g-2$ on the fibers of $C\to T.$ If equality holds, we assume that $D$ is the canonical divisor and $g \geq 2$.

The stability is the (fiber-wise) Gieseker stability for sub-bundles preserved by the Higgs-field, i.e. $(\mathcal E,\Theta)$ is semi-stable if for any $\mathcal F\subseteq\mathcal E$ such that $\Theta\rvert_{\mathcal F}\colon\mathcal F\to\mathcal F\otimes\mathcal O_C(D)$ 
\[\frac{\chi}{r(\mathcal F)}\leq\frac{\chi}{r}.\]
 
 The Hitchin fibration is the morphism:
  \[H\colon \M_{r,\chi}(C)\to B_T:=\bigoplus_{i=1}^r H^0(C,\mathcal O_C(id))\]
 associating to $(\mathcal E,\Theta)$ the coefficients of its characteristic polynomial.

Via the spectral correspondence \cite{hitchin1987stable}, $\BM_{r,\chi}(C)$  can be interpreted as a moduli stack of 1-dimensional sheaves $\mathcal F$ on the non-compact surface \[S_C:=\operatorname{Tot}(\mathcal O_C(D))\xrightarrow{\pi} C\]
with $[\operatorname{Supp}(\mathcal F)]=r[C]$. The hypothesis on the degree of $D$ implies that $S_C$ is del Pezzo, which in turn guarantees the smoothness of the moduli stack $\BM_{r,\chi}(C)$  and its $\BG_m$-rigidification $\FM_{r,\chi}(C)$  for each $(r,\chi)$ similar to Lemma \ref{lem:relativesmooth}.

Moreover, via the correspondence we can identify the fibers of the Hitchin fibration $H$ with the Simpson compactified Jacobians of the spectral curve $C_{P(t)}\subseteq S_C$ cut out by the polynomial $\text{det}(\pi^*\Theta-t\text{Id}).$

In the following theorem we collect the analogous statements of Theorem \ref{thm:semistablesheaves} and Proposition \ref{prop:codim} for the moduli stack of meromorphic Higgs-Bundles. Since the proof is essentially the same as in the previous section we don't repeat it here. 

 


\begin{theorem}\label{summHiggs}

 \begin{enumerate}[leftmargin=0.5cm]
 \item  The moduli stack  $\BM_{r,\chi}(C)$ and its $\mathbb G_m$-rigidification $\FM_{r,\chi}(C)$ are algebraic stacks of finite type over $T$.
 
 \item If $\operatorname{deg}(D) >2g-2$ then  $\BM_{r,\chi}(C)$ is smooth over $T.$
 \item If $D=K_C$ and $T=\operatorname{Spec}(\mathbb{C})$ then  $\BM_{r,\chi}(C)$  is normal.
 
\item There exists a quasi-projective $T$-scheme $\M_{r,\chi}(C)$, proper over $B_T$ which uniformly represent the moduli functor of S-equivalence classes of Higgs bundles.

\item  $\BM^{st}_{r,\chi}(C)$ is a $\BG_m$-gerbe over the moduli space $\M^{st}_{r,\chi}(C)$

\item  $\BM^{st}_{r,\chi}(C)$ is universally closed over $T$ and thus over $\M^{st}_{r,\chi}(C)$

\item The open $\mathcal{U}_{r,\chi}\subseteq\FM_{r,\chi}(C)$  of Higgs bundles whose associated spectral sheaf is a line bundle is dense and its complement has codimension at least two in $\FM_{r,\chi}(C)$. 

\end{enumerate}
\end{theorem}

\begin{rmk}
In the case of Higgs bundle, no extra assumption is necessary to show normality of the moduli stack. We can understand that from at least two points of view: in the language of \cite{davison2022bps}, the moduli stack of Higgs bundles satisfies the \emph{totally negative CY2} property \cite[\S~7]{davison2022bps}; under this additional property the normality follows from the results of \cite{vernet2022rational}. 

Else, we notice that the locus $V$ corresponding to non-reduced spectral curves in the Hitchin base $V$ has always codimension at least two, and thus, as in Proposition~\ref{prop:codim}, the locus of Higgs bundle whose corresponding spectral sheaf is a line bundle has codimension at least two in the stack, which allows us to argue as in \ref{thm:K3case}.

The estimate on the codimension on $V$ can be given by an explicit calculation:

By Riemann-Roch the base $B$ of the Hitchin fibration has dimension 
\[\dim B=d\frac{r(r+1)}{2}+r(1-g),\]
 where $r$ is the fixed rank and $d=\text{deg}\mathcal O_C(D)\geq 2g-2$. The non-reduced locus is the (closure of) the strata  $V_{kr_1,r_2} \subseteq B$ whose points correspond to polynomials in $t$ admitting a  factorisation of the form
\begin{equation}\label{eq:degree}
P(t)=P_1(t)^kP_2(t),\;\; \text{deg}(P_1(t))=r_1\geq 1,\; k\geq 2,\;  \text{deg}(P_2(t))=r_2\geq  0,\;\;\; kr_1+r_2=r.
\end{equation}
Thus the codimension of the non-reduced locus is the minimum of the codimension of the strata $V_{kr_1,r_2} .$
Now, the dimension of the strata is
\[\dim V_{kr_1,r_2}=d\frac{r_1(r_1+1)}{2}+r_1(1-g)+d\frac{r_2(r_2+1)}{2}+r_2(1-g),\]
and it is easy to see by direct computation that
\[\dim B-\dim \dim V_{kr_1,r_2}\geq \frac{4g+1}{2},\]
for $d\geq 2g-1$ and 
\[\dim B-\dim \dim V_{kr_1,r_2}= 3(g-1),\]
for $d=2g-2.$

\end{rmk}

\section{$p$-adic integration}\label{pbasics}
Let $F$ be a non-archimedean local field $F$ with ring of integers $\FO$ and residue field $k \cong \BF_q$.
We consider the normalised Haar measure $\mu$ on the locally compact group $(F,+)$ so that $\mu(\FO) = 1$.

Since $F$ is a completely valued field there is a basic theory of differential geometry available. In particular one can define analytic manifolds and differential forms as over the reals, see \cite{MR1743467,CNS18} for details. 
Given an $n$-dimensional analytic manifold $M$ and a non-vanishing $n$-form $\omega$ on $M$, one defines a Borel measure $\mu_\omega$ on $M$ by writing in a local chart $U \subset F^n$ of $M$
\[\omega_{|U} = f \cdot dx_1\wedge \dots \wedge dx_n,\] 
and integrating $|f|$ against the Haar measure on $F^n$. 

We are interested in analytic manifolds that arise from algebraic geometry. Concretely, given a smooth, separated and finite type $F$-scheme, or more generally $F$-algebraic space $X$, its set of $F$-points $X(F)$ carries naturally the structure of an analytic manifold by the inverse function theorem.

\subsection{The measure on the moduli spaces}\label{gmsm}

Let $\FM$ be a normal algebraic stack with a finite type morphism to $\Spec(\FO)$ of dimension $n$. We assume that there exists a (quasi)-projective variety $\M$ and a morphism $\pi:\FM \rightarrow M$ satisfying the following conditions: there exists a dense open substack $U' \subset \FM$, smooth over $\FO$, that is stabiliser-free i.e. an algebraic space, and that the complement of $U'$ has codimension at least $2$ in $\FM$; there exists a second dense open $U \subset U'$ such that $\pi$ restricted to $U$ induces an isomorphism onto its image in $M$ which we also denote by $U$. 

The goal of this section is to construct a canonical measure on the analytic manifold
\[ M^\natural = M(\FO) \cap U(F). \]

For this we need the following extra assumption: 

\begin{assumption}\label{exprop} For any $x \in M^\natural = M(\FO) \cap U(F)$ there exists a finite extension $L/F$ with ring of integers $\FO_L$ such that the base change $x_L:\Spec(\FO_L) \to M_L = M \times_{\FO} \FO_L$ admits a lift $x'_L:\Spec(\FO_L) \to \FM_L$.
\end{assumption}

Notice that Assumption \ref{exprop} is in particular satisfied if $\FM \to M$ is universally closed \cite[01KA]{stacks-project}.

We fix a Zariski-open cover $U'= \bigcup_{i\in I} U'_i$ such that on each $U'_i$ the relative canonical bundle $\Omega_{U'/\FO}^{n}$ is trivial and pick a trivialising (i.e. non-vanishing) $n$-form $\omega_i$ on $U'_i$. Let $U_i = U \cap U'_i$. We have a decomposition
 \[M^\natural = \bigcup_i M_i^\natural,\] 
with $M_i^\natural = M(\FO) \cap U_i(F)$ and we claim that the measures $\mu_{\omega_i}$ on $M_i^\natural $ glue to a measure on $M^\natural$ that is independent of the choice of $(\omega_i)_i$.

\begin{proposition}[The canonical measure on $M^\natural$]\label{canmes} For all $i,j \in I$ and any measurable subset $A \subset M_i^\natural \cap M_j^\natural$ we have:
\[ \mu_{\omega_i}\left( A   \right) = \mu_{\omega_j}\left( A  \right). \]
In particular the family $(\mu_{\omega_i})_{i\in I}$ glues to a Borel measure $\mu_{can,\FM}$ on $M^\natural$. Furthermore $\mu_{can,\FM}$ is independent of the choice of $(\omega_i)_{i\in I}$.
\end{proposition}

\begin{proof} Since $\omega_i$ and $\omega_j$ are both non-vanishing sections of the same invertible sheaf over $U'$ we have $\omega_{j|U'_i \cap U'_j} = f \cdot \omega_{i|U'_i \cap U'_j}$ with $f \in H^0(U'_i\cap U'_j, \FO^*_{U'_i \cap U'_j})$. We claim that for every $x \in M_i^\natural \cap M_j^\natural \subset  (U_i \cap U_j)(F)$ we have that $f(x) \in \FO_F^*$. This implies the first part of the proposition, since $\mu_{\omega_i}$ and $\mu_{\omega_j}$ are given by locally integrating the absolute value of $\omega_i$ and $\omega_j$. 

To prove the claim, let $L/F$ be a finite extension such that $x_L$ lifts to an $\FO_L$-point $x'_L$ of $\FM_L$ as in Assumption \ref{exprop}. Let $p: X \to \FM$ be a smooth atlas. By passing to an even bigger finite extension of $F$ if necessary, we may further assume that the closed point $\Spec(k_L) \to \FM_{k_L}$ lifts to the atlas $X_{k_L}$. Since $X_L \to \FM_L$ is smooth we can then lift the whole $x'_L$ to a morphism $\tilde{x}_L: \Spec(\FO_L) \to X_L$ by Hensel's Lemma. Since $F^* \cap \FO_L = \FO^*$ it is enough to show that $f(x_L) = f \circ \pi \circ p (\tilde{x}_L) \in \FO_L^*$. 

Write $V_i =  p^{-1}(U'_i)$. The point is now that $f \circ p_{|V_i \cap V_j}$ is a regular non-vanishing function on $V_i \cap V_j$ and the complement of $V_i\cap V_j$ in $X$ has codimension at least $2$ by assumption. Since $X$ is normal, Hartogs' Theorem applies and $f$ thus extends to a regular function on $X$ which is still non-vanishing, as the $0$-locus would have codimension $1$. Therefor
\[ f(x_L) = f \circ p (\tilde{x}_L) \in \FO_L^*,  \]
as claimed. 

The proof of the independence of $\mu_{can,\FM}$ from $(\omega_i)_{i\in I}$ is essentially the same.
\end{proof}

If the stack $\FM$ is clear from the context we will often write $\mu_{can}$ instead of $\mu_{can,\FM}$.

\begin{example} Consider the quotient stack $\FM = [\BA^2/\BG_m]$ where $\BG_m$ acts linearly on $\BA^2$ with weights $(1,-1)$. In the above notation we have
\[ U' = [ (\BA^2 \setminus \{0\}) /\BG_m], \ \ \ U = [(\BA^2 \setminus \{xy=0\}) / \BG_m] \cong \BG_m, \ \ \ M \cong \BA^1. \]
Then $U'$ is isomorphic to the affine line with doubled origin, which admits an open cover by two copies of $\BA^1$. The standard $1$-forms on these $\BA^1$'s glue to a global non-vanishing $1$-form on $U'$ and $\mu_{can}$ is given by integrating $dx$ (i.e. the standard Haar measure) on $M^\natural = \FO \cap F^* = \FO \setminus \{0\}$.
\end{example}

\begin{rmk} If $\FM = M$ is stabiliser-free we have $M^\natural = M(\FO)$ and the construction of $\mu_{can}$ goes back to Weil \cite{weil2012adeles}. The volume is related to the number of $k$-rational points by the formula \cite[Theorem 2.2.5]{weil2012adeles}
\[ \int_{M(\FO)} \mu_{can} = \frac{|M(k)|}{q^n}.\]
\end{rmk}

\begin{rmk} If $\FM$ is Deligne-Mumford stack, then $\pi:\FM \to M$ is proper \cite{MR1432041} and thus Assumption \ref{exprop} is automatically satisfied. In this case $\mu_{can}$ agrees with the orbifold measure constructed in \cite[Section 2.3]{GWZ201}, although in \textit{loc. cit.} there is no assumption on the codimension of the complement of $U \subset \FM$. If $\FM$ is tame, the total volume of $M^\natural$ can be expressed as a weighted point count of the twisted inertia stack of $\FM$ \cite[Theorem 2.21]{GWZ201}.
\end{rmk}

\begin{rmk}\label{rem:measurefromlinebundle} M. Groechenig pointed out to us the following alternative construction: let $\FM/\FO$ be any finite type Artin stack of relative dimension $n$ with a morphism $\pi:\FM \to M$ to a scheme $M$. Assume furthermore, that there exists a smooth stabiliser-free open $U \subset \FM$ such that $\pi_{|U}$ is an isomorphism onto its image and that \ref{exprop} is satisfied. Finally suppose that there exists a line bundle $\widetilde{\mathcal K}$ on $\FM$ which restricts to a power of the canonical line bundle $\mathcal K_{U/\FO}$ on $U$.
Then as in \cite[Section 4.1]{Ya17}, integrating a suitable root of the absolute value of local sections of $\widetilde{\mathcal K}$ against the Haar measure defines a measure $\mu^{vir}$ on $M^\natural$ by an argument as in Proposition~\ref{canmes}.
By construction, if both $\mu^{vir}$ and $\mu_{can}$ are defined they coincide. 

Interestingly, in the context of moduli spaces $\BM$ of sheaves on a surface $S$, such a line bundle $\widetilde{\mathcal K}$ arises from the deformation theory as the determinant of $\bf{R}\mathcal Hom_{\BM}(\mathcal F,\mathcal F)$, where $\mathcal F$ is the universal sheaf on $S\times\BM$ (see \cite{joycedarboux} for the definition in the strictly semi-stable case, or \cite[\S~8.3]{MR1450870} for a simpler explanation on the locus of stables).
\end{rmk}

\subsection{The canonical measure on moduli spaces of sheaves }\label{canmess}

Let $S$ be a smooth projective relative surface over $\Spec(\FO)$ and $\FM_{\beta,\chi}$ the $\BG_m$-rigidified moduli stack of semi-stable 1-dimensional sheaves on $S$ with moduli space $\M_{\beta,\chi}$ as in Section \ref{sec:sheaves}.

We take the open substacks $\mathcal U_{\beta,\chi}^{st}\subset\mathcal U_{\beta,\chi}\subset \FM_{\beta,\chi} $ 
of geometrically stable line bundles and of line bundles respectively as substacks $U\subset U'$  in the notation of Section \ref{gmsm}. From now on we always assume that $S \to \Spec(\FO)$ is either a smooth projective del Pezzo surface or a K3 surface satisfying Assumption \ref{K3as} below:

\begin{assumption}\label{K3as} Let $S \to \Spec(\FO)$ be a smooth projective K3 surface and $(\beta,\chi)$ a dimension vector such that the following holds:
\begin{enumerate} 
\item \label{abcd} The locus of $\mathcal U_{\beta,\chi} \subseteq\mathbb M_{\beta,\chi}(S)$ of line bundles has complement of codimension at least two and $\mathcal U_{\beta,\chi}^{st}$ is non empty.
\item \label{noas} The moduli stack $\mathbb M_{\beta,\chi}(S)$ is normal.
\end{enumerate}
\end{assumption}
\begin{rmk}
We have already commented on the restriction that Assumption \ref{K3as}.\ref{abcd} imposes in Remark \ref{coas}. Over $\BC$, it follows from Theorem \ref{thm:K3case} that Assumption \ref{K3as}.\ref{abcd} implies Assumption \ref{K3as}.\ref{noas}, but we don't know of a proof over more general bases. However, by choosing a suitable spreading out $S \to \Spec(B)$ of a complex K3 surface, satisfying Assumption \ref{K3as}.\ref{abcd}, over the spectrum of a finitely generated $\BZ$-algebra $B$ we obtain finite type moduli stacks $\phi: \BM_{\beta,\chi} \to \Spec(B)$ by Theorem \ref{thm:semistablesheaves}. Then Theorem \ref{thm:K3case} implies that the generic fiber of $\phi$ is geometrically normal and thus by  \cite[Theorem 12.1.6, Corollary 9.5.2]{EGA43}, the fibers of $\phi$ are normal on a non-empty open $V \subset \Spec(B)$. Thus by pulling back $S$ along regular points $\Spec(\FO) \to V$, with $\FO$ the ring of integers of a non-archimedean local field, we obtain many examples of relative K3 surfaces over $\Spec(\FO)$ where Assumption \ref{K3as} holds \cite[Proposition 6.14.1]{EGA42}.

\end{rmk}
By Theorem \ref{thm:semistablesheaves}, Lemma \ref{lem:relativesmooth} and \ref{prop:codim} all the assumptions for Proposition~\ref{canmes} to hold are satisfied and we obtain a canonical measure $\mu_{can}$ on $\M_{\beta,\chi}^\natural$. 

Also the smooth base $B$ of the Hilbert-Chow morphism $h_\chi:\M_{\beta,\chi} \to B$ admits a canonical measure $\mu_{can,B}$ and it will be useful later to have an explicit description of the relative measure $\mu_{can}/\mu_{can,B}$ on smooth fibers of $h_\chi$, as in \cite[Section 6.3] {GWZ20}.

\begin{lemma}\label{relvolf} Let $b \in B(\FO)$ be such that the pullback of the universal curve $\FC_b \to \Spec(\FO)$ is generically smooth. Then under the identification of $h_\chi^{-1}(b)(F)$ with $\Pic^{d}_{\FC_b}(F)$, for $d = \chi + g-1$,  the relative measure $\mu_{can}/\mu_{can,B}$ on $h_\chi^{-1}(b)(F)$ is given by integrating the absolute valued of a translation-invariant gauge form $\omega_b$ on the $\Pic^{0}_{\FC_b}$-torsor $\Pic^{d}_{\FC_b}$. 
\end{lemma}

\begin{proof}
Let $\pi: \operatorname{Pic}^{\text{totdeg}=0}_{\mathcal C /B}\to B$  be 
 the smooth group scheme of total degree zero line bundles on a flat family $\mathcal C\rightarrow B$ of curves. By \cite[\S4.2]{NeronModels} it admits, up to replacing $B$ with a neighborhood of $b\in B$, a global, translation-invariant, trivialising section $\omega_{\text{rel}}$ of $\Omega^g_{\operatorname{Pic}^{\text{totdeg}=0}/B }.$ In particular $\operatorname{Pic}^{\text{totdeg}=0}_{\mathcal C /B}/\mathcal O_F$ has 
a volume form  $\omega=\pi^*\omega_B\wedge\omega_{\text{rel}}$, where $\omega_B$ is a volume form in a neighborhood of $b\in B$ inducing $\mu_{can,B}$.

By  \cite[Lemma~ 6.13] {GWZ20},  since  $\operatorname{Pic}^{\text{totdeg}=d}_{\mathcal C /B}$ is a torsor under the group scheme of total degree zero line bundles,
the relative form $\omega_{\text{rel}}$ induces a section  $\widetilde{\omega}_{\text{rel}}$ of $\Omega^g_{\operatorname{Pic}^{\text{totdeg}=d}/B };$
we thus have a  volume form  $\omega_d=\pi^*\omega_B\wedge\widetilde{\omega}_{\text{rel}}$ on $\operatorname{Pic}^{\text{totdeg}=d}_{\mathcal C /B}$.

Since $\FU_{\beta,\chi}$  has complement of codimension at least two in $\FM_{\beta,\chi}$ by \ref{prop:codim}, integrating the absolute value of $\omega_d$ computes $\mu_{can}$. Therefore on $h_\chi^{-1}(b)(F)$ the relative measure $\mu_{can}/\mu_{can,B}$ is given by integrating the absolute value of $\omega_b = \widetilde{\omega}_{\text{rel},b}$.
\end{proof}

\vspace{0.2cm}

An analogous construction can be done for moduli spaces of Higgs-bundles using Theorem \ref{summHiggs}.




\section{Gerbes and Tate duality}\label{sec:mainsec}

\subsection{$\BG_m$-gerbes}\label{gmgerbes}
Let $X\to V$ be a stack for the fppf-topology. Then $X$ is said to be a \textit{gerbe} if: \begin{enumerate} 
\item for every scheme $T$ over $V$ and objects $x,y \in X(T)$, there exists some fppf-cover $f\colon U \to T$ and an isomorphism $x_{|U}\cong y_{|U}$. 
\item there exists an fppf-cover $U\to V$ such that the groupoid $X(U)$ is non-empty. 
\end{enumerate}
For every $T\to V$ and isomorphism $\phi\colon x \to y$ in $X(T)$, conjugation by $\phi$ induces an isomorphism $\Aut_T(x)\to \Aut_T(y)$. Suppose that every object of $X$ has \textit{abelian} automorphism group in its fibre category; in this case we say that the gerbe $X\to V$ is \textit{abelian}. Then for every $T\to V$ and objects $x,y \in X(T)$, there is a canonical isomorphism $\Aut_T(x)\cong \Aut_T(y)$. These isomorphisms determine a sheaf of groups $Band(X/V)$ over $V$. 

We say that an abelian gerbe $X \rightarrow V$ is a $\BG_m$-gerbe if there is an isomorphism $\BG_{m,V}\to Band(X/V).$ 

By \cite{giraud} equivalence classes of $\BG_m$-gerbes on $V$ are in bijection with the group $H^2(V_{fppf},\BG_m)$, the \textit{Brauer group} of $V$. A gerbe $X/V$ is \textit{trivial} if it is equivalent to $[X/\BG_{m,V}]$, where $\BG_{m,V}$ acts trivially on $X$.

We recall that for $F$ a local field, the Brauer group $H^2(F_{fppf},\BG_m)$ is isomorphic to $\mathbb Q/\mathbb Z$ by mean of the Hasse invariant, see for example \cite[Proposition XIII.6]{MR0354618}.

Given $X/V$ a $\BG_m$-gerbe,  for any point  $F\xrightarrow{s} V$ we get a $\BG_m$-gerbe $s^*X/F$  by pullback. This defines a complex valued function:
\begin{equation}\label{eq:gerbefunction}
\varphi_V\colon V(F)\to H^2(F_{fppf},\BG_m)\cong\mathbb Q/\mathbb Z\xrightarrow{e^{2\pi i} }\mathbb C
\end{equation}

\subsection{Picard stacks/schemes of curves}
Let $V$ be a scheme and $c\colon X\to V$ a smooth proper curve \footnote{i.e. a smooth proper morphism flat and of finite presentation whose fibres are 1-dimensional and connected} with fibres of genus $g\geq 0$.

We denote by $\PIC_{X/V}$ the Picard stack of $X/V$; in other words, for an $V$-scheme $T$, $\PIC_{X/V}(T)$ is the groupoid of invertible sheaves on $X\times_VT$. The group of isomorphism classes of objects of $\PIC_{X/V}(T)$ is therefore the classical Picard group $\Pic(X\times_VT)$. For each integer $d$, we denote by $\PIC^d_{X/V}$ the substack of line bundles of degree $d$. 

We denote by $\Pic_{X/V}$ the sheaf $R^1_{fppf}c_*\BG_m$ on $V$; it can be equivalently defined as the sheafification of the presheaf $T \mapsto \Pic(X\times_VT)$. The natural morphism \begin{equation}\label{eqn:gerbe_pic_pic} \alpha \colon \PIC_{X/V} \to \Pic_{X/V}
\end{equation} is a gerbe \cite[0DME]{stacks-project} banded by $\BG_m$. 

The degree-map factors though $\alpha$ and we denote by $\Pic^d_{X/V}$ the corresponding component of degree $d$ line bundles. The degree zero part $\Pic^0_{X/V}$ is represented by an abelian scheme of relative dimension $g$, the \textit{relative Jacobian}, and each $\Pic^d_{X/V}$ is an $V$-torsor under it.  

In general, a $T$-point of $\Pic_{X/V}$ need not correspond to a line bundle on $X\times_VT$; this property gets lost when the sheafification of $T\mapsto \Pic(X\times_VT)$ is taken.

\begin{lemma}\cite[Proposition 8.4]{NeronModels}\label{lemma:BLR}
For every $V$-scheme $T$, there is a canonical exact sequence of abelian groups, functorial in $T$,
\[0\to \Pic(T) \to \Pic(X\times_VT) \to \Pic_{X/V}(T) \to H^2(T,\BG_m) \to H^2(X\times_VT,\BG_m)\]
\end{lemma}  
Since in our case $\Pic_{X/V}$ is a $V$-scheme, we may look at the particular case $T=\Pic_{X/V}$. Then the identity in $\Pic_{X/V}(\Pic_{X/V})$ gets mapped to the element of $H^2(\Pic_{X/V},\BG_m)$ corresponding to the $\BG_m$-gerbe $\alpha$ of \eqref{eqn:gerbe_pic_pic}.  
By Lemma \ref{lemma:BLR}, $\alpha$ is trivial if and only if there is a universal Poincar\'e bundle $\mathcal L$ on $X\times_V\Pic_{X/V}$. For example, this is the case when $X\to V$ has a section, since then the map \[H^2(\Pic_{X/V},\BG_m) \to H^2(X\times_V\Pic_{X/V},\BG_m)\] is injective.

More generally, for an $V$-scheme $T$, the map $c_T\colon \Pic_{X/V}(T) \to H^2(T,\BG_m)$ sends $x\colon T\to \Pic_{X/V}$ to the equivalence class of the $\BG_m$-gerbe on $T$ obtained by pulling back $\alpha$ along $x$. Such a gerbe is trivial if and only if $x$ comes from a line bundle on $X\times_VT$ (that is, is in the essential image of $\alpha$).

\subsection{Tate duality}\label{tated}
We consider now the special case where the base $V$ is the spectrum of a non-archimedean local field $F$. As recalled above,  Brauer group $H^2(F,\BG_m)$ is identified with $\BQ/\BZ$ via the Hasse invariant. We denote by $f\colon \Pic_{X/F}(F)\to \BQ/\BZ$ the group homomorphism of \ref{lemma:BLR}, and for each integer $d\in \BZ$, we call $f_d\colon \Pic^d_{X/F}(F)\to \BQ/\BZ$ the restriction of $f$. 

Notice that any $x \in \Pic^d_{X/F}(F)$ gives rise to an isomorphism $s_x: \Pic^0_{X/F} \to \Pic^d_{X/F}$ and since $f$ is a group homomorphism we have 
\begin{equation}\label{picgrp}
f_d \circ s_x = f_d(x) + f_0.
\end{equation}

The following proposition gives some control over $f_0$ in terms of torsor instead of gerbes.

\begin{proposition}\label{lemma:smallest_integers}
The image of $f_0$ is $(\frac{1}{d_0}\BZ)/\BZ$, where $d_0$ is the smallest positive integer for which $\Pic^{d_0}_{X/F}(F)$ is non-empty.\\
In particular $f_0$ is constant if and only if $\Pic^{1}_{X/F}(F)$ is non-empty.
\end{proposition}

Similar statements have been used crucially in \cite{GWZ20,GWZ201} and the following discussion is essentially taken from there.

The proof of Proposition \ref{lemma:smallest_integers} relies on Tate duality over local fields. Let $A$ be an abelian variety over $F$ and $A^t$ its dual. 


\begin{theorem}\cite[Theorem 3.7.8]{milneADT}\label{thm:Tate_duality}
There is a canonical perfect pairing
\[A(F)\otimes H^{1}(F,A^t) \to H^2(F,\BG_m)=\BQ/\BZ.\]
\end{theorem}

The pairing admits the following geometric description, see \cite[Remark 3.11]{GWZ20}. First from \cite[Lemma~3.1]{milneADT} we have an isomorphism 
\begin{equation}\label{h12}  H^{1}(F,A^t) \cong \Ext^2(A,\BG_m).\end{equation}
Now, working in the category of commutative group stacks considered in   \cite{brochardduality},  elements of $ \Ext^2(A,\BG_m)$ can be represented by $\BG_m$-gerbes on $A$ with a group structure. 

Given $x \in \A(F)$, $\tau \in H^{1}(F,A^t),$ and $\alpha_\tau \in \Ext^2(A,\BG_m)$ the image of $\tau$ under \eqref{h12}, the paring of $x$ with $\tau$ from Theorem \ref{thm:Tate_duality} equals the class of the $\BG_m$-gerbe $x^*\alpha_\tau$ on $\Spec(F)$ under the Hasse invariant isomorphism $H^2(F,\BG_m) \cong \BQ/\BZ$. 

Thus for $A= \Pic^0_{X/F} = A^t$ the homomorphism $f_0$ is simply the Tate-duality pairing with the $\BG_m$-Gerbe $\alpha_0: \PIC^0_{X/F} \rightarrow \Pic^0_{X/F}$.

\begin{proof}[Proof of Proposition \ref{lemma:smallest_integers}] 
To simplify the notation we drop the subscript $X/F$. By the previous discussion and Theorem \ref{thm:Tate_duality} we see that $d_0\cdot f_0 \equiv 0$ if and only if $\alpha_0^{d_0}$ is equivalent to the trivial $\BG_m$-gerbe on $\Pic^0$. We thus need to show that the latter holds if and only if $\Pic^{d_0}(F)$ is non-empty.

We consider the category of dualizable commutative group stacks $\DCGS_F$ over $\Spec(F)$, see for example \cite{brochardduality} for a detailed account. The internal hom-functor $D(\cdot) =  \sheafhom(\cdot ,B\G_m)$ induces an anti-equivalence on $\DCGS_F$ satisfying $D\circ D = Id$ and extending the usual duality functor on abelian varieties. Furthermore there is an auto-equivalence $D(\PIC) \cong \PIC$ \cite[Section 3.2]{Travkin:2011fk}. Thus applying $D$ to the short exact sequence 
\[0 \to B\G_m \to \PIC \to \Pic \to 0 \]
we get the sequence
\[0 \to D(\Pic) \to \PIC \xrightarrow{\beta} \BZ \to 0. \]
Here exactness on the right follows from \cite[Proposition 3.18]{brochardduality} and the fact that $\mathcal{E}xt^2(\Pic)=0$ \cite[Corollary 11.5]{brochardduality}.
Since $\beta$ is an epimorphism with kernel  $\PIC^0$  we get an equivalence $D(\Pic) \cong \PIC^0$. From this we deduce that the sequences
\begin{align}
\label{osec} 0 \to \Pic^0 \to \Pic \to \BZ \to 0 \\
\label {tsec} 0 \to B\BG_m \to \PIC^0 \to \Pic^0 \to 0
\end{align}
are exchanged by duality.

Now \eqref{osec} is the extension of $\BZ$ by $\Pic^0$ corresponding to the $\Pic^0$-torsor $\Pic^1$ under the equivalence \cite[Proposition 5.8]{brochardduality}. In particular $\Pic^1$ has a rational point if and only \eqref{osec} splits. Dually \eqref{tsec} is the $B\BG_m$-torsor on $\Pic^0$ corresponding to $\alpha_0$ under the equivalence  \cite[Proposition 5.11]{brochardduality} and splits if and only if $\alpha_0$ is trivial. Since  $D(\cdot)$ is an anti-equivalence, we have thus shown that $\Pic^1$ has a rational point if and only if $\alpha_0$ is trivial.

To extend this to $\Pic^{d_0}$ and $\alpha_0^{d_0}$, we notice that their associated short exact sequences are given by the $d_0$-fold Baer sum \cite[010I]{stacks-project} of \eqref{osec} and \eqref{tsec} respectively. By definition $D(\cdot)$ commutes with direct sums and exchanges the addition map with the diagonal, i.e. for every $G \in \DCGS_F$
\[ D\left( G \oplus G \xrightarrow{+} G   \right) \cong D(G) \xrightarrow{\Delta} D(G) \oplus D(G),  \]
and hence also commutes with Baer sums. We thus deduce that $\Pic^{d_0}$ has a rational point if and only if $\alpha_0^{d_0}$ is trivial, which proves the proposition.
\end{proof}

\section{Main theorem} \label{mnthm} We are now ready to state and prove our main theorem. Let  $S \to \Spec(\FO)$ be either a smooth projective del Pezzo surface or a K3 surface satisfying Assumption \ref{K3as} and consider the analytic manifold $\M_{\beta,\chi}^\natural$ together with its canonical measure $\mu=\mu_{can}$ as constructed in Section \ref{canmess}.

The restriction of the natural $\BG_m$-gerbe $\alpha\colon\BM^{st}_{\beta,\chi} \to \FM^{st}_{\beta,\chi}\cong\M^{st}_{\beta,\chi}$ induces a function 
\[ \varphi_{\beta,\chi}: \M^{\natural}_{\beta,\chi}(F) \rightarrow H^2(F,\BG_m) \cong \BQ/\BZ \xrightarrow{e^{2\pi i \cdot}} \BC, \]
by pullback to the generic fiber as defined in \eqref{eq:gerbefunction}.

Now for any $x \in \M_{\beta,\chi}(k)$ denote by $\M_{\beta,\chi}(\FO)_x \subset \M_{\beta,\chi}(\FO)^\natural$ the ball of $\FO$-rational points specializing to $x$ over $k$. Then define the non-archimedean BPS-function $\pBPS_{\beta,\chi}:\M_{\beta,\chi}(k) \rightarrow \BC$ to be:
\begin{equation}\label{eq:pbps}
\pBPS_{\beta,\chi}(x) =q^{-\dim M_{\beta,\chi}}  \int_{M_{\beta,\chi}(\FO)_x} \varphi_{\beta,\chi}^g d\mu_{can}. 
\end{equation}

\begin{theorem}\label{mainthm} The function $\pBPS_{\beta,\chi}: \M_{\beta,\chi}(k) \rightarrow \BC$ satisfies the following two properties:
 \begin{enumerate}
 \item$\pBPS_{\beta,\chi} \equiv q^{-\dim \M_{\beta,\chi}} $ if $(\beta,\chi)$ is generic\footnote{Recall that a pair $(\beta,\chi)$ generic (with respect to $H$), if any Gieseker semi-stable sheaf in $\BM_{\beta,\chi}$ is stable.};
 \item For all $\chi,\chi'\in\mathbb Z$ and for all $y\in B(k)$ we have 
 \[\sum _{x\in h_{\chi}^{-1}(y)(k)} \pBPS_{\beta,\chi}(x)=\sum _{x\in h_{\chi'}^{-1}(y)(k)} \pBPS_{\beta,\chi'}(x).\]
\end{enumerate}
\end{theorem}

Notice that by definition
\[ \sum _{x\in h_{\chi}^{-1}(y)(k)} \pBPS_{\beta,\chi}(x)=\int_{h_\chi^{-1}(B(\FO)_y)} \varphi_{\beta,\chi}^g d\mu.\]
To prove Theorem \ref{mainthm} we want to
 split up the integral using Fubini's theorem. 

First let $B^{sm} \subset B$ be the open subscheme, where the universal curve $\FC \rightarrow B$ is smooth. Then $h_\chi^{-1}( B(\FO)_y\cap B^{sm}(F))$ is contained in $h_\chi^{-1}(B(\FO)_y)$ and its complement has measure $0$, since it is contained in $h_\chi^{-1}((B\setminus B^{sm})(F))$, that is the $F$-points of a closed subscheme of $\M_{\beta,\chi}$ \cite[Proposition 4.4]{GWZ20}. 

Since $h_\chi$ is proper, we can identify the fiber over a point $b \in B(\FO) \cap B^{sm}(F)$ with $h_\chi^{-1}(b)(F)$ where we write $h_\chi^{-1}(b)$ for the $F$-scheme $h_\chi^{-1}(b)_F$. As in \cite[Section 6.3] {GWZ20}  we therefore have an equality:

\begin{equation}\label{fubini}
\sum _{x\in h_{\chi}^{-1}(y)(k)} \pBPS_{\beta,\chi}(x)=\int_{h_\chi^{-1}(B(\FO)_y)} \varphi_{\beta,\chi}^g d\mu = \int_{b\in B(\FO)_y\cap B^{sm}(F)}\left(\int_{h_\chi^{-1}(b)(F)}\varphi_{\beta,\chi}^g d\mu_{\omega_b}\right)d\mu_{can,B}. \end{equation}

Here $\omega_b$ is the gauge form on the $\Pic^0_{\FC_b}$-torsor $h_\chi^{-1}(b) = \Pic^{\chi+g-1}_{\FC_b}$ from Lemma \ref{relvolf}. As the outer integral in \eqref{fubini} is independent of $\chi$ we only need to analyze the fiber integral $\int_{h_\chi^{-1}(b)(F)}\varphi_{\beta,\chi}^g d\mu_{\omega_b}$.

First we notice that the gerbe $\alpha\colon\BM^{st}_{\beta,\chi} \rightarrow \FM^{st}_{\beta,\chi}$ restricted to the fiber over $b$ gets identified with $\PIC^{\chi+g-1}_{\FC_b} \rightarrow \Pic^{\chi+g-1}_{\FC_b}$. Therefore the associated function $\varphi_{\beta,\chi}$ gets identified with 
\[\varphi_{\chi+g-1} = e^{2\pi i f_{\chi+g-1}}\]
where $f_{\chi+g-1}$ is the degree $d=\chi+g-1$ component of the function $f$ introduced in Section \ref{tated}. We may rewrite the inner-most integral in \eqref{fubini} as

\begin{equation}\label{pict} \int_{h_\chi^{-1}(b)(F)}\varphi_{\beta,\chi}^g d\mu_{\omega_b} = \int_{\Pic^{\chi+g-1}_{\FC_b}(F)} \varphi_{\chi+g-1}^g d\mu_{\omega_b}.\end{equation}

Now if $\Pic^{\chi+g-1}_{\FC_b}(F) = \emptyset$ this integral vanishes as we integrate over the empty manifold. Otherwise any $x \in \Pic^{\chi+g-1}_{\FC_b}(F)$ gives an isomorphism  $\Pic^{0}_{\FC_b} \xrightarrow{\sim} \Pic^{\chi+g-1}_{\FC_b}$ and by \eqref{picgrp} we have 
\begin{equation} \label{shift} \int_{h^{-1}(b)(F)}\varphi_{\beta,\chi}^g d\mu_{\omega_b} = \varphi_{\chi+g-1}^g(x) \int_{\Pic^0_{\FC_b}(F)}\varphi_{0}^g d\mu_{\omega_b} \end{equation}

To complete the proof  we need to further study the function $f_0: \Pic^0_{\FC_b}(F) \rightarrow \BQ/\BZ$.

\begin{lemma}\label{lemma:equivalent_conditions}
The following are equivalent:
\begin{enumerate}
\item $f_0 \equiv 0$;
\item $\Pic^1_{\FC_b}(F)$ is non-empty;
\item $\Pic^g_{\FC_b}(F)$ is non-empty.
\end{enumerate}
In this case the image of $f: \Pic_{\FC_b}(F) \rightarrow \BQ/\BZ$ is contained in $(\frac{1}{g} \BZ)/\BZ$.
\end{lemma}

\begin{proof}
The equivalence $(1) \iff (2)$ follows from Proposition \ref{lemma:smallest_integers}. Clearly $(2)$ implies $(3)$ so we only need to prove that $(3)$ implies $(1)$. Suppose then that $\Pic^g_{\FC_b}(F)$ is non-empty. We claim that $f(\Pic^g_{\FC_b}(F))=0$, where 
\[f:\Pic_{\FC_b}(F) \rightarrow \BQ/\BZ\]
 is as is Section \ref{tated}. Consider the moduli stack $\BM_{\beta,1}$ of semi-stable sheaves with Euler characteristic $\chi=1$. For this choice of $\chi$ there are no strictly semi-stable sheaves; hence every semi-stable sheaf is geometrically stable and in particular simple. It follows that the coarse space $M_{\beta,1}$ coincides with the $\BG_m$-rigidification $\FM_{\beta,1}$, hence the coarse moduli map $\BM_{\beta,1}\to \M_{\beta,1}$ gives a $\BG_m$-gerbe $\alpha$ over $\Spec \mathcal O$. For every $x\in \M_{\beta,1}(\mathcal O)$, the pullback $x^*\alpha$ is a trivial gerbe, since $Br(\mathcal O)=0$. In particular, the map 
\[f_{\beta,1}: \M_{\beta,1}(\mathcal O)\to \BQ/\BZ\]
vanishes.
Now notice that line bundles of degree $g$ have Euler charachteristic $1$, i.e. $\Pic^g_{\FC_b}(F)\subset \M_{\beta,1}(\FO)$. This proves the claim. To conclude the proof, let $x\in \Pic^g_{\FC_b}(F)$. For every $y\in Pic^0_{\FC_b}(F)$, $f(y)=f(x)+f(y)=f(x+y)\in f(\Pic^g(F))=0$. Hence $f_0$ is zero.

The statement about the image of $f$ follows from the observation $f(\Pic_{\FC_b}^g(F))=0$.
\end{proof}

\begin{rmk} The equivalence of (2) and (3) in Lemma \ref{lemma:equivalent_conditions} is true more generally. In fact for any geometrically integral, smooth projective curve $C$ of genus $g$ over a local field $\Pic_C^{g-1}$ has a rational point by \cite{Li69} and \cite[Corollary 4]{PS99}.
\end{rmk}

\begin{proof}[Proof of Theorem \ref{mainthm}]
By the previous discussion, in particular \eqref{fubini} and \eqref{pict}, $\chi$-independence amounts to show for all $d,d' \in \BZ$ and $b \in B(\FO) \cap B^{sm}(F)$ the equality

\[\int_{\Pic^{d}_{\FC_b}(F)}\varphi_{d}^g d\mu_{\omega_b} = \int_{\Pic^{d'}_{\FC_b}(F)}\varphi_{d'}^g d\mu_{\omega_b}.\]

For simplicity take $d' = g$. Consider first the case where $\Pic^g_{\FC_b}(F)=\emptyset$. Then the right-hand side is zero. If $\Pic^d_{\FC_b}(F)$ is also empty, we are done. Else let $x\in \Pic^d(F)$. Then as in \eqref{shift}
\begin{equation}\label{shiftt} \int_{\Pic^{d}_{\FC_b}(F)}\varphi_{d}^g d\mu_{\omega_b} = \varphi_{d}^g(x) \int_{\Pic^0_{\FC_b}(F)}\varphi_{0}^g d\mu_{\omega_b}.\end{equation}

As $\Pic_{\FC_b}^g(F)$ is empty, so is $\Pic_{\FC_b}^1(F)$ by Lemma \ref{lemma:equivalent_conditions}, hence $f_0$ is a non-trivial character, i.e. a surjective homomorphism $\Pic_{\FC_b}^0(F)\to \BZ/d_0\BZ$ for some $d_0\geq 2$ by Lemma \ref{lemma:smallest_integers}. Furthermore $d_0$ is the smallest integer such that $\Pic_{\FC_b}^{d_0}(F)\neq \emptyset$, hence $d_0$ does not divide $g$. The function $\varphi_0^g$ is therefore non-zero and by a character-sum argument the right hand side of \eqref{shiftt} vanishes, as required.

Next we consider the case $\Pic_{\FC_b}^g(F)\neq \emptyset$. Then by Lemma \ref{lemma:equivalent_conditions} also $\Pic_{\FC_b}^1(F)\neq \emptyset$ and thus also $\Pic_{\FC_b}^d(F)\neq \emptyset$. Furthermore the image of $f$ is $g$-torsion and thus $\varphi_d^g=\varphi_g^g \equiv 1$. Since $\Pic_{\FC_b}^d(F) \cong \Pic_{\FC_b}^0(F) \cong \Pic_{\FC_b}^g(F)$ we are done with the proof of $(2)$.

The proof of (1) follows from the fact that for $(\beta,\chi)$ generic, the stable locus is everything, thus the morphism $\BM_{\beta,\chi} \to \M_{\beta,\chi}$ is a $\BG_m$-gerbe. In this case $f_{\beta,\chi} \equiv 0$, since the gerbe $x^*\alpha \in H^2(F,\BG_m)$ is pulled back from $H^2(\Spec(\FO),\BG_m) = 0$ and thus trivial.
\end{proof}

\begin{rmk} It is straightforward to extend the definition of the $\pBPS$-function and the proof of Theorem \ref{mainthm} to the case of usual and meromorphic Higgs bundles using Theorem \ref{summHiggs}.
\end{rmk}

\bibliographystyle{amsalpha}

\bibliography{master}
\end{document}